\documentclass[reqno]{amsart}

\usepackage{float}
\usepackage[utf8]{inputenc}

\usepackage{caption}
\usepackage{capt-of}
\usepackage[usenames,dvipsnames,svgnames,table]{xcolor}
\usepackage[colorlinks=true, pdfstartview=FitV, linkcolor=blue, citecolor=blue, urlcolor=darkblue]{hyperref}

\usepackage{geometry}                
\usepackage{graphicx}
\usepackage{mathrsfs}
\usepackage{amssymb}
\usepackage{amsmath}
\usepackage{bbm}
\usepackage{amsfonts}
\usepackage{mathrsfs}
\usepackage{epstopdf}
\usepackage{lscape}
\usepackage{array}
\usepackage[utf8]{inputenc}
\usepackage{tikz,caption}
\usepackage[maxbibnames=99]{biblatex}
\usepackage{enumitem}
\setlist[itemize]{leftmargin=2em}
\setlist[enumerate]{leftmargin=2em}

\usepackage{booktabs}
\usepackage{multirow}
\usepackage{mathtools}
\usepackage[linesnumbered,ruled]{algorithm2e}
\usepackage{tikz}
\usepackage{mathtools}
\usepackage{soul}
\usepackage{upgreek}

\newtheorem{theorem}{Theorem}[section]
\newtheorem{lemma}[theorem]{Lemma}
\newtheorem{corollary}[theorem]{Corollary}

\newtheorem{proposition}[theorem]{Proposition}

\theoremstyle{definition}
\newtheorem{definition}[theorem]{Definition}
\newtheorem{example}[theorem]{Example}

\newtheorem{remark}[theorem]{Remark}

\newcolumntype{C}[1]{>{\centering\arraybackslash}m{#1}}

\newcommand{\Z}{\mathbb{Z}}
\newcommand{\C}{\mathbb{C}}

\newcommand{\Q}{\mathbb{Q}}

\newcommand{\pc}{\mathsf{PC}}
\newcommand{\csf}{\mathbf{X}}

\newcommand{\csff}{\mathbf{X}_F}
\newcommand{\csft}{\mathbf{X}_T}
\newcommand{\stfrak}{\mathfrak{st}}
\newcommand{\lmulti }{\{\!\!\{}
\newcommand{\rmulti}{\}\!\!\}}

\newcommand{\sort}{\mathrm{sort}}
\newcommand{\lc}{\lambda_\mathrm{LC}}
\newcommand{\lead}{\lambda_\mathrm{lead}}

\definecolor{darkblue}{rgb}{0.0,0,0.7} 
\definecolor{darkred}{rgb}{0.7,0,0} 
\definecolor{darkgreen}{rgb}{0, .6, 0} 

\newcommand{\defncolor}{\color{darkred}}
\newcommand{\defn}[1]{{\defncolor\emph{#1}}} 

\usepackage[colorinlistoftodos]{todonotes}

\newcommand{\tikzmath}[2][]{\vcenter{\hbox{\begin{tikzpicture}[#1]#2
\end{tikzpicture}}}
}

\title{Proper implies improper for the chromatic symmetric function}

\author[R.~Orellana]{Rosa Orellana}
\address[R.~Orellana]{Mathematics Department, Dartmouth College,
Hanover, NH, USA}
\email{Rosa.C.Orellana@dartmouth.edu}
\urladdr{\href{https://math.dartmouth.edu/~orellana/}{https://math.dartmouth.edu/~orellana/}}

\author[M.~Tomba]{Mario Tomba}
\address[M.~Tomba]{Mathematics Department, University of Chicago,
Chicago, IL, USA}
\email{mariotomba@uchicago.edu}
\urladdr{\href{https://mariotomba.github.io/}{https://mariotomba.github.io/}}

\date{}

\begin{document}

\begin{abstract}
    Stanley's question of whether the chromatic symmetric function distinguishes trees has been unanswered for over thirty years. It has now been checked for trees with at most 29 vertices, and for several classes of trees. There is the perception that answering the question might be easier for the class of ``proper'' trees, i.e. trees where every vertex is either a leaf or adjacent to a leaf.  In this paper, we show that an affirmative answer for a subclass of proper trees, which we call ``$k$-weighted'', is equivalent to an affirmative answer for all trees.   In addition, we show that in the class of $k$-weighted trees, the coefficients carry redundant information. In fact we give a closed formula for obtaining the redundant coefficients from other coefficients. 
\end{abstract}

\maketitle

\section{Introduction}
The chromatic symmetric function is a generalization of the chromatic polynomial of a graph that retains substantially more information about the graph's structure. Introduced by Stanley in 1995 \cite{stanley95symmetric}, it has become an important object at the intersection of algebraic combinatorics, graph theory, and the theory of symmetric functions. Let $G=(V,E)$ be a finite simple graph with  $V=\{v_1, \ldots, v_n\}$. A proper coloring of $G$ is a map $\kappa: V\rightarrow \mathbb{P}$ such that adjacent vertices receive different values. The \emph{chromatic symmetric function} of $G$ is defined by 
 \[ \csf_G(\boldsymbol{x})= \sum_{\kappa \in \pc(G)} x_{\kappa(v_1)}\cdots x_{\kappa(v_n)}\]
where $\boldsymbol{x}= (x_1, x_2, \ldots)$ is an infinite set of commuting variables and $\pc(G)$ is the set of proper colorings of $G$.  

The chromatic symmetric function generalizes the chromatic polynomial of a graph, since setting $x_1= \cdots = x_k =1$ and all other variables to zero gives the chromatic polynomial $\chi_G(k)$, which counts the number of proper colorings of $G$ using at most $k$ colors. In general, two non-isomorphic graphs can have the same chromatic polynomial but different chromatic symmetric functions. Thus, $\csf_G$ is a stronger invariant than $\chi_G(k)$. However, there are known examples of pairs of non-isomorphic graphs with equal chromatic symmetric function, which raises the natural question of classifying those classes of graphs that are distinguished by their chromatic symmetric function. 

One central open problem in the area is the \emph{Tree Isomorphism Conjecture}, which asserts that if $T$ and $T'$ are finite trees satisfying $\csf_T = \csf_{T'}$, then $T$ and $T'$ must be isomorphic. In other words, the chromatic symmetric function is a complete isomorphism invariant for trees. The question of whether the chromatic symmetric function distinguishes trees was first raised in 1995 by Stanley in \cite{stanley95symmetric} and more formally in \cite[Problem 7.47(b)]{stanley1999enumerative}. 

Since being proposed, the Tree Isomorphism Conjecture has drawn a lot of attention. Computationally, Heil and Ji \cite{heil2019algorithm} verified it for all trees up to 29 vertices. More generally, it has been shown to hold for several subclasses of trees that satisfy certain structural properties. A \emph{proper} tree is a tree where every vertex is either a leaf or adjacent to a leaf, if a tree is not proper, we call it \emph{improper}. A \emph{spider} is a tree with at most one vertex of degree greater than two, while a \emph{2-spider} is a related class with two branching vertices. A \emph{caterpillar} is a tree that becomes a path after all leaves are removed.  The conjecture has been shown to hold for the following classes of trees: proper caterpillars \cite{aliste2014proper}, spiders \cite{martin2008distinguishing, HurynChmutov20}, all caterpillars \cite{loebl2018isomorphism}, trees with exactly two vertices of degree at least 3 \cite{wang2024class}, proper trees with diameter at most 5 \cite{ADOZ}, and for all trees with diameter at most 5 \cite{GOT26}, trees that are determined by their combinatorial Laplacian spectra \cite{AWvW24}, proper $q$-caterpillars \cite{ANRS24}, and proper trees whose weighted core is multiplicity-isolated \cite{Awan26}.

The caterpillar case was extended in stages. In 2008, Martin, Morin and Wagner proved the conjecture for a subclass of proper caterpillars \cite{martin2008distinguishing}. In 2014, Aliste-Prieto and Zamora \cite{aliste2014proper} extended the result to all proper caterpillars, and later in 2019, Loebl and Sereni \cite{loebl2018isomorphism} settled the case of all caterpillars by relating it to an isomorphism problem for weighted trees. Similarly, for trees of diameter at most five, the case of proper trees was proved first in \cite{ADOZ}, while the general case was settled later by the authors of this paper in \cite{GOT26} using the deletion-near-contraction relation. 

There are structural reasons for showing the case of proper trees first. We say an edge is \emph{internal} if neither endpoint is a leaf, and that a non-leaf vertex is \emph{deep} if it is not adjacent to any leaf; a tree is proper exactly when it has no deep vertices. For an internal edge $e$, the deletion-near-contraction (DNC) relation, see Section \ref{subsection: DNC} for details, gives
a recursive algorithm that terminates when we reach a forest with no internal edges (a \emph{star forest}), and this yields the star-basis expansion of $\csf_T$. If $T$ is proper, each of the three operations in the DNC relation removes exactly one internal edge and the resulting forests are again proper. The computation results in a complete ternary tree of depth $m=|I(T)|$, the number of internal edges of $T$. Moreover, since no cancellations occur \cite{ADOZ}, the absolute values of the coefficients sum to $3^m$. However, when $T$ is improper, a single deletion can remove two or even three internal edges. In this case, the computation of $\csft$ depends on the adjacencies between deep vertices in $T$, and the coefficients are much harder to control. For instance, this difference is already visible in the leading term. For every proper tree $T$, the smallest partition $\lambda$ in lexicographic order occurring in the star expansion of $\csft$ always has coefficient $c_\lambda = 1$.  However, when $T$ has $d$ deep vertices, the leading coefficient depends on the degrees of these vertices and the leading partition has $d$ parts of size 1.

In fact, our main result shows that if the Tree Isomorphism Conjecture is true for proper trees, then it will be true for all trees. Indeed, we have further reduced the Tree Isomorphism Conjecture to trees where every non-leaf vertex is adjacent to exactly $k-1$ leaves, for any fixed $k\geq 2$. We call this class of trees \emph{$k$-weighted trees}. The key ingredient is Theorem \ref{thm:general attachment csf}, which expresses the chromatic symmetric function of a $k$-weighted tree as an evaluation of the chromatic symmetric function of its core.

Moreover, we show that in the case of $k$-weighted trees, we only need a subset of the coefficients to verify the Tree Isomorphism conjecture. In Theorem \ref{thm:1-free coefficients}, we show that the coefficients indexed by partitions containing ones can be written in terms of those coefficients that do not contain ones. Hence, if the Tree Isomorphism conjecture is true for $k$-weighted trees, we will be able to determine it from the coefficients indexed by 1-free partitions. 

This paper is organized as follows. In Section \ref{sec:preliminaries}, we review basic definitions and results in graph theory and symmetric functions. In Section \ref{sec:csf and star basis}, we review the deletion-near-contraction relation, how it allows the computation of the chromatic symmetric functions in terms of the \emph{core} of a graph, and the structural information revealed by the leading partition. In Section \ref{sec:proper tree reduction}, we review the rooted product operation on graphs, and use it to prove our main result, that is, that verifying the conjecture for proper trees implies the full conjecture. Lastly, in Section \ref{sec:k weighted trees}, we prove that the 1-free coefficients of $k$-weighted trees determine their entire chromatic symmetric function.

\subsection*{Acknowledgments}  R. Orellana was partially supported by NSF grant DMS-2452044.

\section{Preliminaries and Definitions}\label{sec:preliminaries}

This section contains the background necessary in graph theory and symmetric functions for the results in this paper. We will follow West's book \cite{westBook} in graph theory and Stanley's Enumerative Combinatorics 2 \cite[Chapter 7]{stanley1999enumerative} for symmetric functions. 

\subsection{Graph theory background}
 A \defn{graph} $G$ is a pair $(V,E)$, where $V$ is the set of \defn{vertices} and $E$ is a multiset of \defn{edges}, where an \defn{edge} $e$ is an unordered pair of vertices $uv$, in which case $u$ and $v$ are called the \defn{endpoints} of $e$. In the case that we are discussing several graphs, we will use $V(G)$ and $E(G)$ to indicate that we are referring to the vertex and edge set of $G$, respectively.

 The \defn{order} of a graph is the number of vertices, i.e. $|V|$. We say that an edge $e$ is \defn{incident} to a vertex $v$ if $v$ is an endpoint of $e$. The \defn{degree} of a vertex $v$, denoted by $\deg(v)$, is the number of edges incident to it. An edge $uu$ is called a \defn{loop}. We say that $G$ is \defn{simple} when it contains no loops and there can be at most one edge between two distinct vertices.  In other words, $G$ is simple when $E$ is a \emph{set} containing no loops.
    
    Two graphs $G, H$ are \defn{isomorphic} if there exists a bijection $\varphi:V(G) \to V(H)$ between the vertex sets of the graphs such that $uv \in E(G)$ if and only if $\varphi(u)\varphi(v) \in E(H)$. We say that $G$ is \defn{connected} if for any $u,v \in V(G)$, there is a sequence of distinct vertices $v_1,v_2, \ldots,v_n$ where $u=v_1$ and $v=v_n$, such that $v_iv_{i+1}\in E(G)$. A \defn{cycle} in $G$ is a sequence of distinct edges in which only the first and last vertex are equal.

    A \defn{tree} $T$ is an acyclic connected graph.  A vertex of degree 1 in a tree is called a \defn{leaf}. We say that a tree is \defn{proper} if every vertex is either a leaf or adjacent to at least one leaf.  A \defn{forest} is a graph where each connected component is a tree.

    Let $\mathbb{P}$ denote the positive integers. A \defn{proper coloring} of $G$ is a function $\kappa : V(G) \to \mathbb{P}$ such that $\kappa(u) \neq \kappa(v)$ for any $uv \in E$. We let $\pc(G)$ denote the set of proper colorings of $G$. When the image $\kappa(V)$ has size $|\kappa(V)| \leq k$, we say that $\kappa$ is a \defn{$k$-coloring}. The number of $k$-colorings is denoted by $\chi_G(k)$. It is a well-known result in graph theory that the function $\chi_G : V \to \mathbb{P}$ is a polynomial of degree $|V|$, called the \defn{chromatic polynomial}.

    We now introduce some definitions around special types of edges and vertices that we will use repeatedly. An \defn{internal vertex} is a vertex of degree at least two.  An\defn{internal edge} is an edge incident to two internal vertices, and a \defn{leaf edge} is an edge such that one of its endpoints is a leaf. We denote by $I(G)$ the set of internal edges of $G$ and by $IV(G)$ the set of internal vertices of $G$.
    \begin{figure}[ht]
        \centering
        \begin{tikzpicture}[auto=center,every node/.style={circle, fill=black, scale=0.5}]
            \node (a) at (1,0) {};
            \node (b) at (0,.5) {};
            \node (c) at (0,-.5) {};
            \node (d) at (-1,.5) {};
            \node (e) at (-1,-.5) {};
            \node (f) at (2, .5) {};
            \node (g) at (2, -.5) {};
            \node (h) at (3, .5) {};
            \node (i) at (3, -.5) {};

            \draw[thick] (a)--(b);
            \draw[thick] (a)--(c);
            \draw[thick] (b)--(c);
            \draw[thick] (b)--(d);
            \draw[thick] (c)--(e);
            \draw[thick] (f)--(h);
            \draw[thick] (g)--(a)--(f)--(g);
            \draw[thick] (g)--(i);

            \node[fill=none, scale=1.75] at (-.5, .7) {$\ell$};
            \node[fill=none, scale=1.75] at (-1.2, .6) {$v$};
            \node[fill=none, scale=1.75] at (1, 0.2) {$u$};
            \node[fill=none, scale=1.75] at (2.2, 0) {$e$};
        \end{tikzpicture}
        \caption{The edge $\ell$ is a leaf edge, while $e$ is an internal edge. The vertex $v$ is a leaf vertex, and $u$ is an internal vertex.}
        \label{fig:basic graph concepts}
    \end{figure}
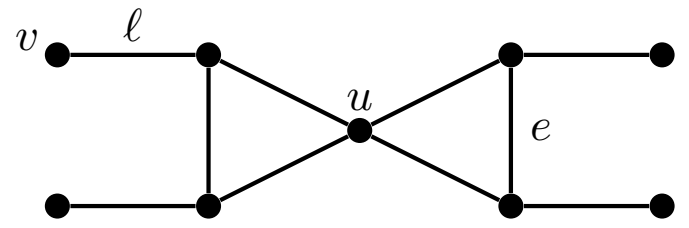

\subsection{Symmetric functions}\label{subsec:sym-fun} We now recall some definitions about symmetric functions. We begin with partitions, since they index bases of symmetric functions. 

    A \defn{partition} $\lambda$ is a sequence of positive integers $(\lambda_1, \ldots, \lambda_k)$ such that $\lambda_1 \geq \cdots \geq \lambda_k$. 
    The quantity $|\lambda| \coloneqq \sum_i \lambda_i$ is called the \defn{weight} of $\lambda$. If $|\lambda|=n$, we say that $\lambda$ is a partition of $n$, and denote it $\lambda \vdash n$. Each $\lambda_i$ is a \defn{part}, and the number of parts is the \defn{length} of $\lambda$, denoted $\ell(\lambda)$. There is a natural total ordering for partitions called \defn{lexicographic order}, where whenever $\lambda, \mu \vdash n$, we say that $\lambda \leq \mu$ if $\lambda = \mu$ or if $\lambda_i = \mu_i$ for $1\leq i<j$ and $\lambda_j<\mu_j$ for some $1\leq j \leq \ell(\lambda)$. If the partitions have different lengths, we set $\lambda_i = 0$ if $i>\ell(\lambda)$, similarly for $\mu$.

    Let $\boldsymbol{x}= (x_1, x_2, \ldots)$ be a countably infinite set of commuting variables. We say $f \in \Q[[x_1,x_2,\ldots]]$ is a \defn{symmetric function} if $f$ is invariant under the action of the symmetric group $S_\infty =\bigsqcup_{n \geq 1} S_n$, which acts by permuting the indices of the variables and $f$ has bounded degree.
    
    For $r\geq 1$, let $p_r$ denote the \defn{$r$-th power sum symmetric function} which is defined by 
    \[ p_r \coloneqq x_1^r + x_2^r + x_3^r+ \cdots \in \Q[[x_1,x_2,\ldots]] \]
    and for a partition $\lambda = (\lambda_1, \ldots, \lambda_\ell)$, we define the \defn{power sum symmetric function} by 
    \[ p_\lambda \coloneqq p_{\lambda_1} \cdots p_{\lambda_\ell}~.\]
    It is well-known that the set $\{p_\lambda \, |\, \lambda\vdash k\}$ is linearly independent and that $p_\lambda$ is symmetric for any partition $\lambda$. Then, the \defn{algebra of symmetric functions} is a graded subalgebra of $\Q[[x_1, x_2, \ldots]]$ defined as
    \[ 
    \Lambda \coloneqq \bigoplus_{n \geq 0}\Lambda^n ~,
    \]
    where $\Lambda^0 \coloneqq\mathbb{Q}$ and $\Lambda^k \coloneqq \text{Span}_{\mathbb{Q}} \{ p_\lambda \, |\, \lambda \vdash k\}$ for $k \ge 1$.

\section{The Chromatic Symmetric Function and the star-basis}\label{sec:csf and star basis}
    In this paper, we are interested in  the \defn{chromatic symmetric function} introduced by Stanley in \cite{stanley95symmetric}. For any finite simple graph $G$ with vertex set $\{v_1, \ldots, v_n\}$ define
    \[
    \csf_G = \sum_{\kappa \in \pc(G)} x_{\kappa(v_1)}x_{\kappa(v_2)}\cdots x_{\kappa(v_n)}.
    \]
    Note that $\csf_G$ is a homogeneous symmetric function of degree $n = |V(G)|$. In \cite{stanley95symmetric}, Stanley described the expansion of $\csf_G$ in terms of the power sum basis. Given $A \subseteq E(G)$, the \defn{spanning subgraph} induced by $A$ is the subgraph $G|_A$ with vertex set $V$ and edge set $A$. If $G$ has $n$ vertices, let $\lambda(A)$ be the partition of $n$ whose parts are the orders of the connected components of $G|_A$.

    \begin{theorem}[\cite{stanley95symmetric}, Theorem 2.5]\label{thm:stanley-power-sum-expansion}
        We have
        \[
        \csf_G = \sum_{A \subseteq E} (-1)^{|A|}p_{\lambda(A)}
        \]
    \end{theorem}
    
    In \cite{chromBases}, Cho and van Willigenburg introduced a plethora of symmetric function bases using the chromatic symmetric function. In this paper, we only need a special case of the general result that was proved in that paper. For the general statement, we refer the reader to \cite[Theorem 5]{chromBases}.

    For any positive integer $k$, let $St_k$ denote the tree with $k-1$ vertices of degree 1 and one vertex of degree $k-1$. This graph is called the \defn{star graph} on $k$ vertices or \defn{$k$-star}. For example, for $k=5$, we have
    \[
        St_5 = 
            \tikzmath{
            \filldraw (0,0) node {} circle (.06cm);
            \filldraw (0,-.6) node {} circle (.06cm);
            \filldraw (0,.6) node {} circle (.06cm);
            \filldraw (-.6,0) node {} circle (.06cm);
            \filldraw (.6,0) node {} circle (.06cm);
        
            \draw[thick] (0,0)--(0,-.6);
            \draw[thick] (0,0)--(0,.6);
            \draw[thick] (0,0)--(0.6,0);
            \draw[thick] (0,0)--(-0.6,0);
    }
    \]

    For a partition $\lambda = (\lambda_1, \ldots, \lambda_\ell)$, define $St_\lambda = St_{\lambda_1} \sqcup \cdots \sqcup St_{\lambda_\ell}$ as the disjoint union of $\lambda_i$-stars. Then, define 
    \[
    \stfrak_k \coloneqq \csf_{St_k} \text{ and }  \stfrak_\lambda \coloneqq \stfrak_{\lambda_1}\cdots\stfrak_{\lambda_\ell}~.
    \] 

    \begin{theorem}[\cite{chromBases}] \label{thm:cho-vanW}
        For any integer $n \geq 1$, the set $\{\stfrak_\lambda : \lambda \vdash n\}$ is a basis for $\Lambda^n$.
    \end{theorem}

We call $\{\stfrak_\lambda : \lambda \vdash n\}$ the \defn{star-basis}, which has been recently studied in the context of Stanley's tree isomorphism conjecture in \cite{ADOZ, GOTfpscac,GOT26}. Moreover, there is a nice change of basis from the power sum basis to the star-basis that was proved in \cite{chromReciprocity}.

\begin{proposition} 
    We have
\begin{equation}
        \stfrak_{n+1} = \sum_{r=0}^n(-1)^r \binom{n}{r} p_{(r+1,1^{n-r})}
    \end{equation}
    and 
    \begin{equation}\label{eq:power-sum-to-star}
        p_{n+1}=\sum_{r=0}^n(-1)^r \binom{n}{r} \stfrak_{(r+1,1^{n-r})}
    \end{equation}
\end{proposition}

\subsection{Deletion-Near-Contraction} \label{subsection: DNC} In this section, we review the \defn{deletion-near-contraction} (DNC) relation for the chromatic symmetric function, which was proved in \cite{ADOZ}, and further utilized in \cite{GOTfpscac, GOT26}. The DNC relation is a recursive relation that allows for an efficient computation of the chromatic symmetric function in the star-basis. For an internal edge $e\in E(G)$, the three graph-theoretic operations are the following:

\begin{itemize}
    \item {\bf Deletion:} this is the usual deletion of an edge in a graph. For an edge $e \in E(G)$, we denote the graph obtained by deleting $e$ by $G \setminus e$.
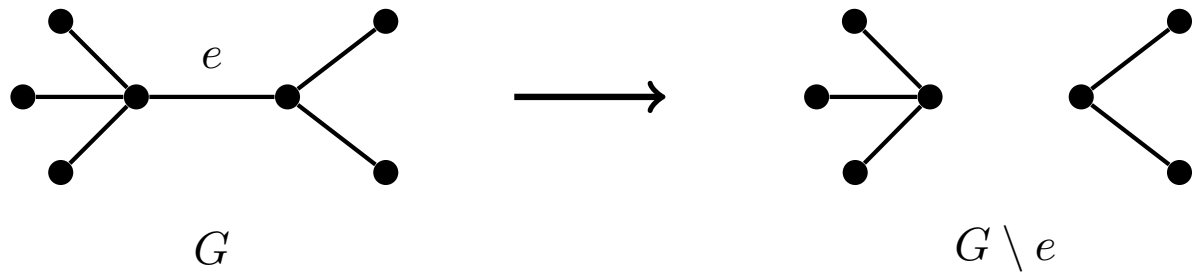
    \begin{center}
        \begin{tikzpicture}[auto=center,every node/.style={circle, fill=black, scale=0.5}]
            \node (a) at (0,0) {};
            \node (a1) at (-.5,.5) {};
            \node (a2) at (-.75, 0) {};
            \node (a3) at (-.5, -.5) {};

            \draw[thick] (a) -- (a1);
            \draw[thick] (a) -- (a2);
            \draw[thick] (a) -- (a3);

            \node (b) at (1,0) {};
            \node (b1) at (1.65, .5) {};
            \node (b2) at (1.65, -.5){};

            \draw[thick] (b) -- (b1);
            \draw[thick] (b) -- (b2);

            \draw[thick] (a) -- node[above, fill=none, scale=1.75] {$e$} ++ (b); 

            \node[fill=none, scale=1.75] at (0.5, -1) {$G$};

            \draw[very thick, ->] (2.5,0) -- (3.5,0);

            \node (a) at (5.25,0) {};
            \node (a1) at (4.75,.5) {};
            \node (a2) at (4.5, 0) {};
            \node (a3) at (4.755, -.5) {};

            \draw[thick] (a) -- (a1);
            \draw[thick] (a) -- (a2);
            \draw[thick] (a) -- (a3);

            \node (b) at (6.25,0) {};
            \node (b1) at (6.9, .5) {};
            \node (b2) at (6.9, -.5){};

            \draw[thick] (b) -- (b1);
            \draw[thick] (b) -- (b2);

            \node[fill=none, scale=1.75] at (5.75, -1) {$G\setminus e$};
        \end{tikzpicture}
    \end{center}
 \item \textbf{Leaf-contraction:} Given $e \in E(G)$, the \defn{leaf-contracted graph} is obtained by contracting $e$ and connecting a leaf with an edge $\ell_e$ to the vertex that results from the contraction of $e$. The leaf-contracted graph is denoted $G\odot e$,
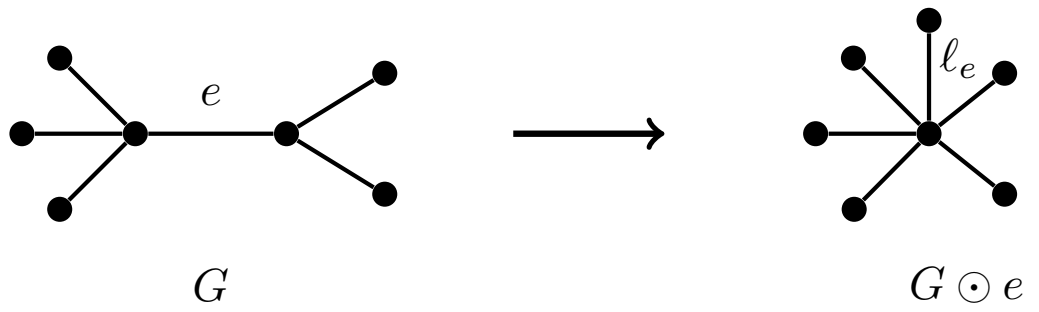
    \begin{center}
        \begin{tikzpicture}[auto=center,every node/.style={circle, fill=black, scale=0.5}]
            \node (a) at (0,0) {};
            \node (a1) at (-.5,.5) {};
            \node (a2) at (-.75, 0) {};
            \node (a3) at (-.5, -.5) {};

            \draw[thick] (a) -- (a1);
            \draw[thick] (a) -- (a2);
            \draw[thick] (a) -- (a3);

            \node (b) at (1,0) {};
            \node (b1) at (1.65, .4) {};
            \node (b2) at (1.65, -.4){};

            \draw[thick] (b) -- (b1);
            \draw[thick] (b) -- (b2);

            \draw[thick] (a) -- node[above, fill=none, scale=1.75] {$e$} ++ (b); 

            \node[fill=none, scale=1.75] at (0.5, -1) {$G$};

            \draw[very thick, ->] (2.5,0) -- (3.5,0);

            \node (a) at (5.25,0) {};
            \node (a1) at (4.75,.5) {};
            \node (a2) at (4.5, 0) {};
            \node (a3) at (4.755, -.5) {};

            \draw[thick] (a) -- (a1);
            \draw[thick] (a) -- (a2);
            \draw[thick] (a) -- (a3);

            \node (b1) at (5.75, .4) {};
            \node (b2) at (5.75, -.4){};

            \node (l) at (5.25, 0.75) {};

            \draw[thick] (a)--(l);

            \draw[thick] (a) -- (b1);
            \draw[thick] (a) -- (b2);

            \node[fill=none, scale=1.75] at (5.45, .5) {$\ell_e$};
            \node[fill=none, scale=1.75] at (5.5, -1) {$G\odot e$};
        \end{tikzpicture}
    \end{center}
 
    \item \textbf{Dot-contraction:} Given $e \in E(G)$, the \defn{dot-contracted graph} is obtained by contracting the edge $e$ and adding an isolated vertex, $v$, to the graph. This can be formulated in terms of the leaf-contraction operation as simply removing the edge $\ell_e$, and thus we denote the dot-contracted graph by $(G\odot e)\setminus \ell_e$,
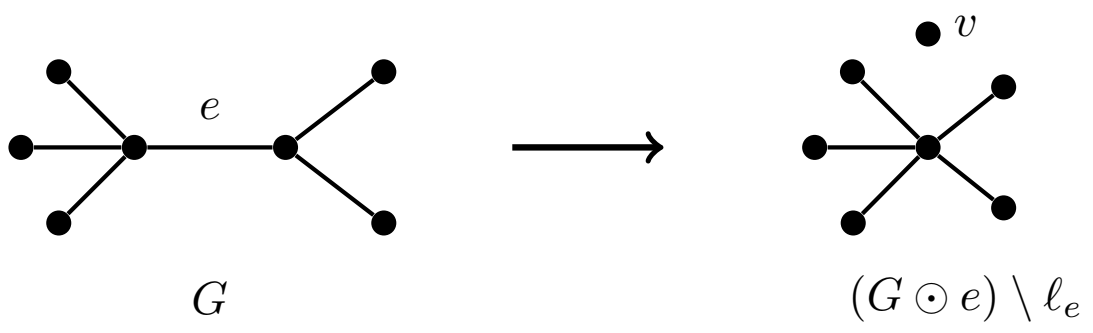
    \begin{center}
        \begin{tikzpicture}[auto=center,every node/.style={circle, fill=black, scale=0.5}]
            \node (a) at (0,0) {};
            \node (a1) at (-.5,.5) {};
            \node (a2) at (-.75, 0) {};
            \node (a3) at (-.5, -.5) {};

            \draw[thick] (a) -- (a1);
            \draw[thick] (a) -- (a2);
            \draw[thick] (a) -- (a3);

            \node (b) at (1,0) {};
            \node (b1) at (1.65, .5) {};
            \node (b2) at (1.65, -.5){};

            \draw[thick] (b) -- (b1);
            \draw[thick] (b) -- (b2);

            \draw[thick] (a) -- node[above, fill=none, scale=1.75] {$e$} ++ (b); 

            \node[fill=none, scale=1.75] at (0.5, -1) {$G$};

            \draw[very thick, ->] (2.5,0) -- (3.5,0);

            \node (a) at (5.25,0) {};
            \node (a1) at (4.75,.5) {};
            \node (a2) at (4.5, 0) {};
            \node (a3) at (4.755, -.5) {};

            \draw[thick] (a) -- (a1);
            \draw[thick] (a) -- (a2);
            \draw[thick] (a) -- (a3);

            \node (b1) at (5.75, .4) {};
            \node (b2) at (5.75, -.4){};

            \node at (5.25, 0.75) {};

            \draw[thick] (a) -- (b1);
            \draw[thick] (a) -- (b2);

            \node[fill=none, scale=1.75] at (5.5, .8) {$v$};
            \node[fill=none, scale=1.75] at (5.5, -1) {$(G\odot e)\setminus \ell_e$};
        \end{tikzpicture}
        \vspace{-.5cm}
    \end{center}
\end{itemize}

The recurrence using these operations is as follows. 

\begin{proposition}[\cite{ADOZ}](The DNC relation) \label{prop:dnc-formula}
    For a finite simple graph $G$ and any edge $e \in E(G)$, we have
    \[
    \csf_G = \csf_{G \setminus e} - \csf_{(G \odot e)\setminus \ell_e} + \csf_{G \odot e}
    \]
\end{proposition}
Each of the three operations defined above reduces the number of internal edges by at least one. 
Since star forests are the only graphs with no internal edges, this relation yields a recursive algorithm to compute the star-basis expansion of the chromatic symmetric function, as proved in \cite{ADOZ}. In addition, applying the relation on a leaf edge $e$ yields $\csf_G = \csf_{G \odot e}$ and $\csf_{G \setminus e} = \csf_{(G \odot e)\setminus \ell_e}$. Hence, we only apply it to internal edges.

   In \cite{ADOZ}, the authors show that we can recursively apply the DNC relation on internal edges until $\csf_G$ can be written as a linear combination of $\csf_H$, where $H$ is a forest of star graphs.  This process is formalized in the \emph{star-expansion algorithm} presented in \cite{ADOZ}. The output of the algorithm is a tree with root $G$ called a \defn{DNC tree} for $G$, \defn{$\mathcal{T}(G)$}. Each leaf of $\mathcal{T}(G)$ is a forest of stars. The orders of the connected components of each forest is a partition that indexes a summand in $\csf_G$. For details about this algorithm see \cite[Section 3]{GOT26}. 

\begin{example}\label{exa:DNC-tree}
    The figure below shows an example of how to apply the star-expansion algorithm. We use red to indicate the internal edge on which we apply the DNC relation. 
    
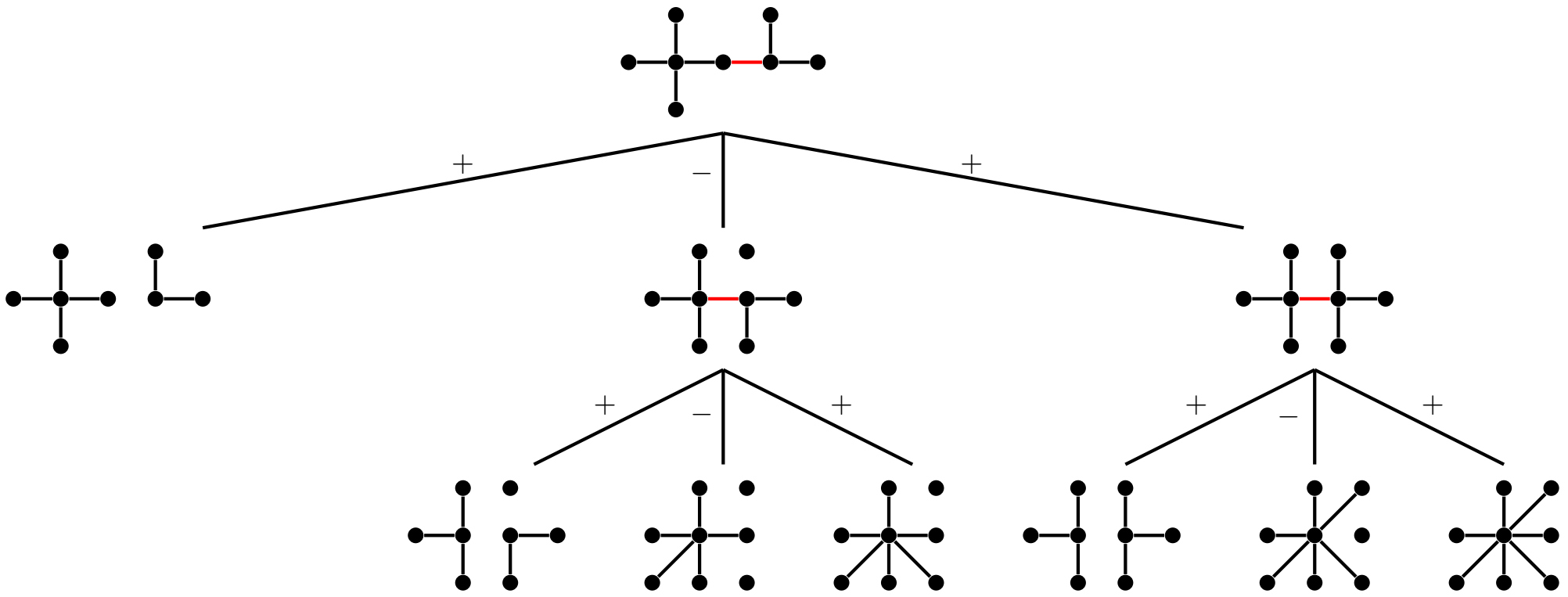
\begin{center}
    \begin{tikzpicture}[auto=center, every node/.style={circle, fill=black, scale=0.4}, thick, scale=0.4]
    \node (A1) at (0,0) {};
    \node (A2) at (1,0) {};
    \node (A3) at (1,1) {};
    \node (A4) at (1,-1) {};
    \node (A5) at (2,0) {};
    \node (A6) at (3,0) {};
    \node (A7) at (4,0) {};
    \node (A8) at (3,1) {};

    \draw (A1)--(A2);
    \draw (A2)--(A3);
    \draw (A2)--(A4);
    \draw (A2)--(A5);
    \draw[red] (A5)--(A6);
    \draw (A6)--(A7);
    \draw (A6)--(A8);

    \node (B1) at (-13,-5) {};
    \node (B2) at (-12,-5) {};
    \node (B3) at (-12,-4) {};
    \node (B4) at (-12,-6) {};
    \node (B5) at (-11,-5) {};
    \node (B6) at (-10,-5) {};
    \node (B7) at (-9,-5) {};
    \node (B8) at (-10,-4) {};

    \draw (B1)--(B2);
    \draw (B2)--(B3);
    \draw (B2)--(B4);
    \draw (B2)--(B5);
    \draw (B6)--(B7);
    \draw (B6)--(B8);

    \node (C1) at (0.5,-5) {};
    \node (C2) at (1.5,-5) {};
    \node (C3) at (1.5,-4) {};
    \node (C4) at (1.5,-6) {};
    \node (C5) at (2.5,-5) {};
    \node (C6) at (3.5,-5) {};
    \node (C7) at (2.5,-4) {};
    \node (C8) at (2.5,-6) {};

    \draw (C1)--(C2);
    \draw (C2)--(C3);
    \draw (C2)--(C4);
    \draw[red] (C2)--(C5);
    \draw (C5)--(C6);
    \draw (C5)--(C8);

    \node (D1) at (13,-5) {};
    \node (D2) at (14,-5) {};
    \node (D3) at (14,-4) {};
    \node (D4) at (14,-6) {};
    \node (D5) at (15,-5) {};
    \node (D6) at (16,-5) {};
    \node (D7) at (15,-4) {};
    \node (D8) at (15,-6) {};

    \draw (D1)--(D2);
    \draw (D2)--(D3);
    \draw (D2)--(D4);
    \draw[red] (D2)--(D5);
    \draw (D5)--(D7);
    \draw (D5)--(D6);
    \draw (D5)--(D8);

    \node (E1) at (-4.5,-10) {};
    \node (E2) at (-3.5,-10) {};
    \node (E3) at (-3.5,-9) {};
    \node (E4) at (-3.5,-11) {};
    \node (E5) at (-2.5,-10) {};
    \node (E6) at (-1.5,-10) {};
    \node (E7) at (-2.5,-9) {};
    \node (E8) at (-2.5,-11) {};

    \draw (E1)--(E2);
    \draw (E2)--(E3);
    \draw (E2)--(E4);
    \draw (E5)--(E6);
    \draw (E5)--(E8);

    \node (F1) at (0.5,-10) {};
    \node (F2) at (1.5,-10) {};
    \node (F3) at (1.5,-9) {};
    \node (F4) at (1.5,-11) {};
    \node (F5) at (2.5,-10) {};
    \node (F6) at (2.5,-11) {};
    \node (F7) at (2.5,-9) {};
    \node (F8) at (0.5,-11) {};

    \draw (F1)--(F2);
    \draw (F2)--(F3);
    \draw (F2)--(F4);
    \draw (F2)--(F5);
    \draw (F2)--(F8);

    \node (G1) at (4.5,-10) {};
    \node (G2) at (5.5,-10) {};
    \node (G3) at (5.5,-9) {};
    \node (G4) at (5.5,-11) {};
    \node (G5) at (6.5,-10) {};
    \node (G6) at (6.5,-11) {};
    \node (G7) at (6.5,-9) {};
    \node (G8) at (4.5,-11) {};

    \draw (G1)--(G2);
    \draw (G2)--(G3);
    \draw (G2)--(G4);
    \draw (G2)--(G5);
    \draw (G2)--(G6);
    \draw (G2)--(G8);

    \node (H1) at (8.5,-10) {};
    \node (H2) at (9.5,-10) {};
    \node (H3) at (9.5,-9) {};
    \node (H4) at (9.5,-11) {};
    \node (H5) at (10.5,-10) {};
    \node (H6) at (11.5,-10) {};
    \node (H7) at (10.5,-9) {};
    \node (H8) at (10.5,-11) {};

    \draw (H1)--(H2);
    \draw (H2)--(H3);
    \draw (H2)--(H4);
    \draw (H5)--(H7);
    \draw (H5)--(H6);
    \draw (H5)--(H8);

    \node (I1) at (13.5,-10) {};
    \node (I2) at (14.5,-10) {};
    \node (I3) at (14.5,-9) {};
    \node (I4) at (14.5,-11) {};
    \node (I5) at (15.5,-10) {};
    \node (I6) at (15.5,-11) {};
    \node (I7) at (15.5,-9) {};
    \node (I8) at (13.5,-11) {};

    \draw (I1)--(I2);
    \draw (I2)--(I3);
    \draw (I2)--(I4);
    \draw (I2)--(I6);
    \draw (I2)--(I7);
    \draw (I2)--(I8);

    \node (J1) at (17.5,-10) {};
    \node (J2) at (18.5,-10) {};
    \node (J3) at (18.5,-9) {};
    \node (J4) at (18.5,-11) {};
    \node (J5) at (19.5,-10) {};
    \node (J6) at (19.5,-11) {};
    \node (J7) at (19.5,-9) {};
    \node (J8) at (17.5,-11) {};

    \draw (J1)--(J2);
    \draw (J2)--(J3);
    \draw (J2)--(J4);
    \draw (J2)--(J5);
    \draw (J2)--(J6);
    \draw (J2)--(J7);
    \draw (J2)--(J8);

    \draw (2,-1.5)--(-9,-3.5);
    \node[fill=none, scale=1.75] at (-3.5,-2.15) {$+$};

    \draw (2,-1.5)--(2,-3.5);
    \node[fill=none, scale=1.75] at (1.55,-2.35) {$-$};

    \draw (2,-1.5)--(13,-3.5);
    \node[fill=none, scale=1.75] at (7.25,-2.15) {$+$};

    \draw (2,-6.5)--(-2,-8.5);
    \node[fill=none, scale=1.75] at (-0.5,-7.25) {$+$};

    \draw (2,-6.5)--(2,-8.5);
    \node[fill=none, scale=1.75] at (1.55,-7.45) {$-$};

    \draw (2,-6.5)--(6,-8.5);
    \node[fill=none, scale=1.75] at (4.5,-7.25) {$+$};

    \draw (14.5,-6.5)--(10.5,-8.5);
    \node[fill=none, scale=1.75] at (12,-7.25) {$+$};

    \draw (14.5,-6.5)--(14.5,-8.5);
    \node[fill=none, scale=1.75] at (13.95,-7.5) {$-$};

    \draw (14.5,-6.5)--(18.5,-8.5);
    \node[fill=none, scale=1.75] at (17,-7.25) {$+$};
    \end{tikzpicture}
    \end{center}

    In particular, it says that for the graph $T$ at the root, we have 
    \[
    \csft = 
    -\mathfrak{st}_{(4,3,1)}+ \mathfrak{st}_{(4,4)}+\mathfrak{st}_{(6,1,1)}+\mathfrak{st}_{(5,3)}-2\mathfrak{st}_{(7,1)}+\mathfrak{st}_{(8)}~.
    \]
\end{example}

If we denote three operations $L$ for deletion, $M$ for dot-contraction, and $R$ for leaf-contraction, then every path in $\mathcal{T}(G)$ from the vertex $G$ to a star forest $F$ can be encoded by a sequence of $L$, $M$, and $R$. More formally, for a star forest $F$, let \defn{$\lambda(F)$} be the partition whose parts are the orders of the connected components of $F$. In this case, we say that $F$ is a forest \defn{indexed by} $\lambda(F)$. Then, the following holds.

\begin{lemma}\label{lem:coeff-paths-seq} Let $G$ be an $n$-vertex graph and $\mathcal{T}(G)$ a DNC-tree corresponding to $G$. If $\csf_G=\sum_{\lambda\vdash n} c_\lambda \mathfrak{st}_\lambda$, then 
\[ c_\lambda = (-1)^{m} |\mathcal{S}_\lambda|,\]
where $\mathcal{S}_\lambda$ is the set of sequences from the root $G$ to a forest $F$ such that $\lambda(F)= \lambda$ and $m$ is the number of $M$s in each sequence.
\end{lemma}

Note that Lemma \ref{lem:coeff-paths-seq} allows us to compute coefficients using combinatorial arguments for the number of paths in a DNC tree, which will be relevant in Section \ref{sec:k weighted trees}.

\subsection{Cores} Recall that for a tree $T$, $IV(T)$ is the set of internal vertices of $T$. For each $v \in V(T)$, let $L_v$ denote the number of leaves in the neighborhood of $v$. As defined in \cite{ADOZ}, the \defn{(weighted) core} of $T$ is the weighted graph $(K_T,w)$ where 
    \begin{itemize}
        \item $K_T=T|_{IV(T)}$ is the subgraph induced by the internal vertices of $T$, and
        \item $w(v) = 1+L_v$ for $v \in V(K)$.
    \end{itemize}
When there is no confusion, we will simply write $K=K_T$ for the unweighted core.

\begin{example}
    A tree $T$ and its weighted core $(K,w)$ are shown below. Note that we are representing the weight of the vertices in the core with vertex labels.
    
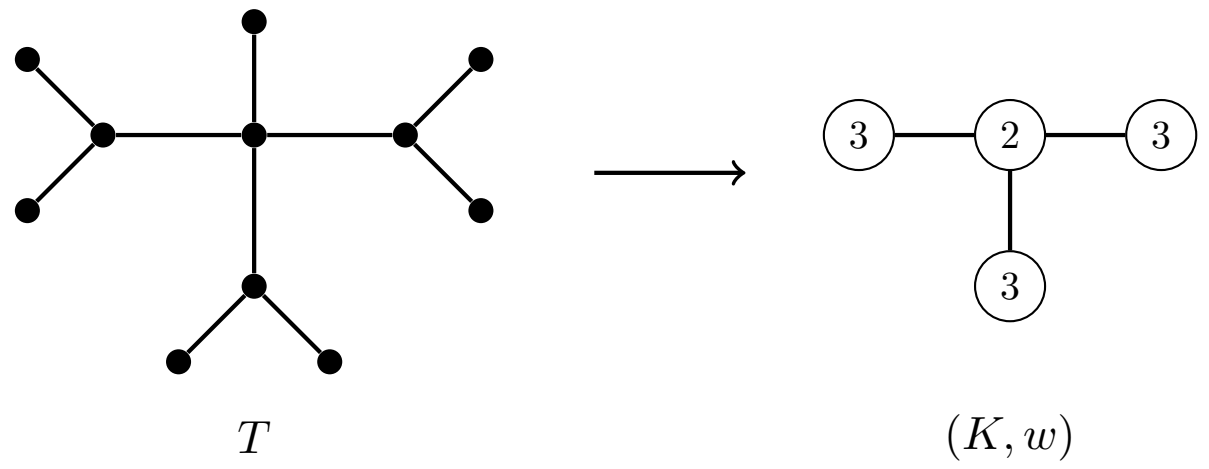
\begin{center}
        
        \begin{tikzpicture}[auto=center,every node/.style={circle, fill=black, scale=0.5}]
            \node (a) at (0,0) {};
            \node (b) at (1,0) {};
            \node (c) at (2,0) {};
            \node (d) at (1,-1) {};
        
            \node (a1) at (-.5,.5) {};
            \node (a2) at (-.5,-.5) {};
            \node (b1) at (1,.75) {};
            \node (c1) at (2.5,.5) {};
            \node (c2) at (2.5,-.5) {};
            \node (d1) at (.5,-1.5) {};
            \node (d2) at (1.5,-1.5) {};
        
            \draw[thick] (a)--(b)--(c);
            \draw[thick] (b)--(d);
            \draw[thick] (a)--(a1);
            \draw[thick] (a)--(a2);
            \draw[thick] (b)--(b1);
            \draw[thick] (c)--(c1);
            \draw[thick] (c)--(c2);
            \draw[thick] (d)--(d1);
            \draw[thick] (d)--(d2);
        
            \node[fill=none, scale=1.75] at (1,-2) {$T$};
        
            \draw[thick, ->] (3.25,-.25)--(4.25,-.25);
        
            \node[fill=white, draw, scale=1.5] (A) at (5,0) {$3$};
            \node[fill=white, draw, scale=1.5] (B) at (6,0) {$2$};
            \node[fill=white, draw, scale=1.5] (C) at (7,0) {$3$};
            \node[fill=white, draw, scale=1.5] (D) at (6,-1) {$3$};
        
            \draw[thick] (A)--(B)--(C);
            \draw[thick] (B)--(D);
        
            \node[fill=none, scale=1.75] at (6,-2) {$(K,w)$};
        \end{tikzpicture}
    \end{center}
\end{example}

In \cite{ADOZ}, the authors proved that the chromatic symmetric function of a tree allows an expansion in terms of the core. For a subset of edges $A \subseteq E(G)$ of a graph, let $c(A)$ denote the set of connected components of $G|_A$. 
In addition, we let $w(B) = \sum_{v \in V(B)}w(v)$.

\begin{proposition}[\cite{ADOZ}, Corollary 5.11]\label{prop:csf in terms of core} Let $T$ be a tree, and $(K, w)$ its weighted core. Then
    \[
    \csft = \sum_{A \subseteq E(K)} \prod_{B \in c(A)} \left( \sum_{j=0}^{|V(B)|-1} (-1)^j \binom{|V(B)|-1}{j} \stfrak_{(w(B)-j,1^j)} \right)
    \] 
    Moreover, if $T$ is proper, this expression is cancellation-free.
\end{proposition}

\begin{example}\label{ex:csf in terms of core}
    We show an example of how to use Proposition \ref{prop:csf in terms of core} to compute $\csft$. Consider the tree $T$ and its core $(K,w)$ below:
    
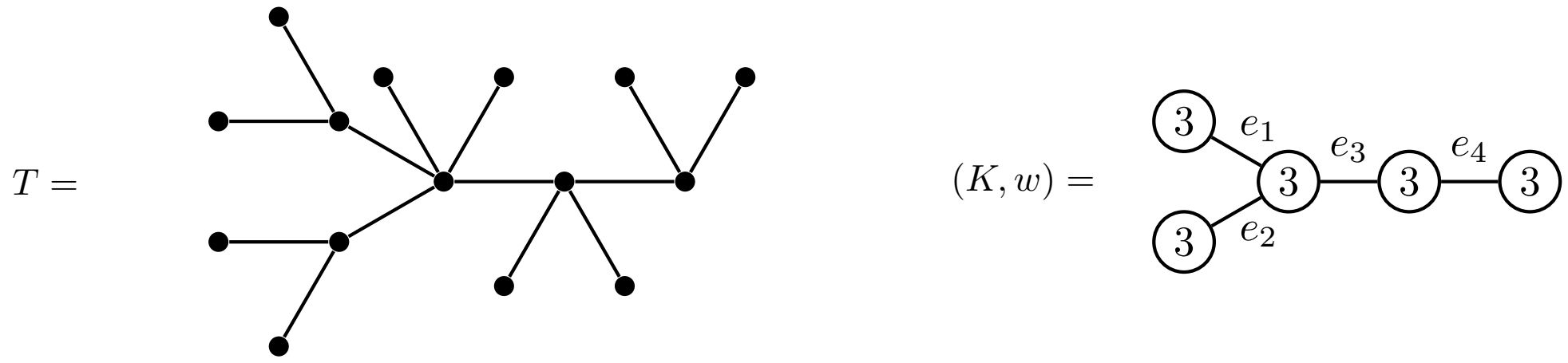
\begin{center}
        \begin{tikzpicture}[
            auto=center,
            every node/.style={circle, fill=black, scale=0.5},
            thick
            ]
            \node (p1) at (0,0) {};
            \node (p2) at (-0.866,0.5) {};
            \node (p3) at (-0.866,-0.5) {};
            \node (p4) at (1,0) {};
            \node (p5) at (2,0) {};
        
            \node (a1) at (-0.5,0.866) {};
            \node (b1) at (0.5,0.866) {};
            \node (a2) at (-1.366,1.366) {};
            \node (b2) at (-1.866,0.5) {};
            \node (a3) at (-1.866,-0.5) {};
            \node (b3) at (-1.366,-1.366) {};
            \node (a4) at (0.5,-0.866) {};
            \node (b4) at (1.5,-0.866) {};
            \node (a5) at (1.5,0.866) {};
            \node (b5) at (2.5,0.866) {};
        
            \draw (p1)--(p2);
            \draw (p1)--(p3);
            \draw (p1)--(p4)--(p5);
        
            \foreach \i in {1,...,5}
                \draw (a\i)--(p\i)--(b\i);
        
            \node[fill=none, scale=1.75] at (-3.3,0) {$T =$};
        
            \begin{scope}[
                xshift=7cm,
                every node/.style={
                    circle,
                    draw,
                    fill=white,
                    minimum size=5mm,
                    inner sep=1pt
                }
            ]
                \node (c1) at (0,0) {$3$};
                \node (c2) at (-0.866,0.5) {$3$};
                \node (c3) at (-0.866,-0.5) {$3$};
                \node (c4) at (1,0) {$3$};
                \node (c5) at (2,0) {$3$};
        
                \draw (c1)-- node[draw=none, fill=none, above right] {$e_1$} (c2);
                \draw (c1)-- node[draw=none, fill=none, below right] {$e_2$} (c3);
                \draw (c1)-- node[draw=none, fill=none, above] {$e_3$} (c4);
                \draw (c4)-- node[draw=none, fill=none, above] {$e_4$} (c5);
        
                \node[draw=none, fill=none, scale=0.875] at (-2.2,0)
                    {$(K,w)=$};
            \end{scope}
        \end{tikzpicture}
    \end{center}
    Given a fixed $A \subseteq E(K)$ and $B \in c(A)$, observe that $w(B) = 3|V(B)|$, since every vertex of $K$ has weight $3$. Then, by Proposition \ref{prop:csf in terms of core},
    \[
    X_T
        =
        \sum_{A\subseteq E(K)}
        \prod_{B\in c(A)}
        \left(
            \sum_{j=0}^{|V(B)|-1}
            (-1)^j\binom{|V(B)|-1}{j}
            \stfrak_{(3|V(B)|-j,1^j)}
        \right).
    \]
    Now note that the expression above only depends on the orders of the connected components in $K|_A$ for a given $A \subseteq E(K)$. For example, the subsets $A_1= \{e_1\}$ and $A_2=\{e_2\}$ both yield the same term in the summation since $\lambda(A_1)=\lambda(A_2) = (2,1^3)$. In particular, each $A_1$ and $A_2$ contribute $
    \stfrak_{(6,3^3)} - \stfrak_{(5,3^3,1)}$ to $\csft$.
    
    Thus, we group the edge subsets according to $\lambda(A)$
    \[
    \begin{array}{c|c|c}
    |A| & \lambda(A) & \text{Number of such } As \\ \hline
    0 & (1^5)   & 1 \\
    1 & (2,1^3) & 4 \\
    2 & (3,1^2) & 4 \\
    2 & (2^2,1) & 2 \\
    3 & (4,1)   & 3 \\
    3 & (3,2)   & 1 \\
    4 & (5)     & 1
    \end{array}.
    \]
    Then,
    \[
\begin{aligned}
X_T ={}&
\stfrak_{(3^5)}
+4\left(\stfrak_{(6)}-\stfrak_{(5,1)}\right)
    \stfrak_{(3^3)} +4\left(
    \stfrak_{(9)}-2\stfrak_{(8,1)}+\stfrak_{(7,1^2)}
\right)\stfrak_{(3^2)} \\
&+2\left(\stfrak_{(6)}-\stfrak_{(5,1)}\right)^2
    \stfrak_{(3)} +3\left(
    \stfrak_{(12)}-3\stfrak_{(11,1)}
    +3\stfrak_{(10,1^2)}-\stfrak_{(9,1^3)}
\right)\stfrak_{(3)} \\
&+\left(
    \stfrak_{(9)}-2\stfrak_{(8,1)}+\stfrak_{(7,1^2)}
\right)
\left(\stfrak_{(6)}-\stfrak_{(5,1)}\right) \\
&+\stfrak_{(15)}-4\stfrak_{(14,1)}
+6\stfrak_{(13,1^2)}-4\stfrak_{(12,1^3)}
+\stfrak_{(11,1^4)} 
\end{aligned}
\]
\end{example}

\subsection{The leading partition} Here we review some of the results obtained in \cite{GOT26} regarding the smallest partition in lexicographic order that appears with nonzero coefficient in the chromatic symmetric function. Let $F$ be an $n$-vertex forest with chromatic symmetric function $\csf_F=\sum_{\lambda \vdash n}c_\lambda \mathfrak{st}_\lambda$. The \defn{leading partition} of $\csf_F$ is the smallest partition $\lambda \vdash n$, in lexicographic order, such that $c_{\lambda} \neq 0$. We then say that $c_{\lambda}$ is the \defn{leading coefficient} and that $c_{\lambda} \mathfrak{st}_{\lambda}$ is the \defn{leading term}. We denote the leading partition of $F$ by $\lead(\csff)$.

Recall that $I(F)$ is the set of internal edges of $F$. We refer to the connected components of $F \setminus I(F)$ as the \defn{leaf components} of $F$, and we denote by $\lc(F)$ the partition whose parts are the orders of the leaf components of $F \setminus I(F)$:
\[
\lc(F)\coloneqq\lambda(F\setminus I(F)).
\]
We will call this partition the \defn{leaf component partition} of $F$.
Notice that $F \setminus I(F)$ is a spanning subgraph, therefore every vertex in $F$ is also a vertex in $F \setminus I(F)$. Furthermore, $F\setminus I(F)$ is a forest whose connected components are all stars. Therefore, a leaf component is always a star tree. 

\begin{theorem}[\cite{GOT26}, Theorem 4.15]\label{thm:leading partition} For any forest $F$, we have
\[
\lc(F) = \lead(\csff)
\]   
\end{theorem}

\begin{remark}\label{rmk:lead partition and core weights}
    Given a tree $T$ and its core $(K,w)$, note that for any $v \in V(K)$, $w(v)$ is the order of the leaf component containing $v$ in $T$. Thus, by Theorem \ref{thm:leading partition}, the parts of $\lead(\csft)$ are precisely the weights $\lmulti w(v) : v\in V(K)\rmulti$.
\end{remark}

\begin{example}
    Let $T$ be the tree from Example \ref{ex:csf in terms of core}. Then, simplifying the star-basis expansion for $\csft$, we obtain:
    \[
        \begin{aligned}
            X_T ={}&
            \stfrak_{(3^5)}
            -4\stfrak_{(5,3^3,1)}
            +2\stfrak_{(5^2,3,1^2)}
            +4\stfrak_{(6,3^3)}
            -4\stfrak_{(6,5,3,1)}
            +2\stfrak_{(6^2,3)} \\
            &+4\stfrak_{(7,3^2,1^2)}
            -\stfrak_{(7,5,1^3)}
            +\stfrak_{(7,6,1^2)}
            -8\stfrak_{(8,3^2,1)}
            +2\stfrak_{(8,5,1^2)}
            -2\stfrak_{(8,6,1)} \\
            &-3\stfrak_{(9,3,1^3)}
            +4\stfrak_{(9,3^2)}
            -\stfrak_{(9,5,1)}
            +\stfrak_{(9,6)}
            +9\stfrak_{(10,3,1^2)}
            +\stfrak_{(11,1^4)} \\
            &-9\stfrak_{(11,3,1)}
            -4\stfrak_{(12,1^3)}
            +3\stfrak_{(12,3)}
            +6\stfrak_{(13,1^2)}
            -4\stfrak_{(14,1)}
            +\stfrak_{(15)}.
        \end{aligned}
    \]
    Then, $\lead(\csft)=(3^5)=\lc(T)$, which also coincides with $\lmulti w(v) : v\in K\rmulti$.
\end{example}

In \cite{GOT26}, the authors show several properties about the leading partition, and it is a key ingredient in the reconstructions they obtain for trees of small-diameter.
Here we only need one of those results. Recall that given a partition $\lambda$, $m_p(\lambda)$ is the multiplicity of $p$ in $\lambda$.
\begin{corollary}[\cite{GOT26}, Corollary 4.20]
    A tree $T$ is proper if and only if $\lead(\csft)$ has no parts of size 1, i.e., $m_1(\lead(\csft)) = 0$.
\end{corollary}

\section{Proper-tree reduction}\label{sec:proper tree reduction}

In this section, we prove our main result, i.e., we reduce the Tree Isomorphism conjecture to a subclass proper trees. We will begin by reviewing the \emph{rooted product of graphs}, an operation on graphs that first appeared in the work by Godsil and McKay \cite{GM} in the context of graph spectra, and then we will prove the main theorem. 

\subsection{The rooted product of graphs} We review the rooted product of graphs, denoted $G(H)$, and show that $G\cong G'$ if and only if $G(H_r)\cong G'(H_r)$. Recall that a \defn{rooted graph} is a pair $(H,r)$ where $H$ is a simple graph and $r \in V(H)$ is a distinguished vertex, called the \defn{root}.

\begin{definition}\label{def: graph attachment}
    Let $G$ be a finite simple graph, and $(H,r)$ a rooted connected graph. We define their \defn{rooted product}, denoted by $G(H_r)$ to be the graph with vertex set
    \[
    V(G(H_r)) = V(G) \times V(H),
    \]
    and edges given by
    \begin{equation}\label{eq:attachment edge set}
    E(G(H_r)) = \{(g,r)(g',r) : gg' \in E(G)\} \cup \{(g,h)(g,h') : hh'\in E(H), g \in V(G)\}.
    \end{equation}
    Informally, $G(H_r)$ is obtained by attaching a copy of $H$ at each vertex $v \in V(G)$ and identifying the root $r$ with $v$. When the root $r$ is clear from the context, we simply denote the rooted product by $G(H)$.
\end{definition}

\begin{example}\label{ex:graph attachment}
    We provide an example for the operation defined above. Let $G$ and $(H,r)$ be the following graphs:
    
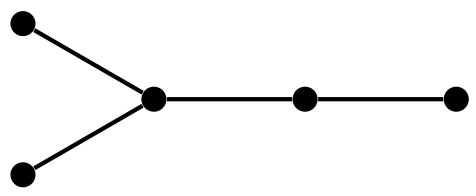
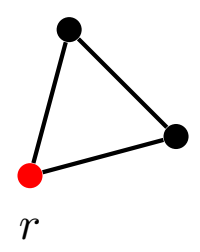
\begin{center}
        \begin{tikzpicture}[auto=center,every node/.style={circle, fill=black, scale=0.5}]
            \node (g1) at (0,0) {};          
            \node (g2) at (-0.866,0.5) {};   
            \node (g3) at (-0.866,-0.5) {};  
            \node (g4) at (1,0) {};          
            \node (g5) at (2,0) {};          
        
            \draw[thick] (g1)--(g2);
            \draw[thick] (g1)--(g3);
            \draw[thick] (g1)--(g4)--(g5);
        
            \node[fill=none, scale=1.75] at (-1.9,0) {$G =$};
        
            \node[fill=none, scale=1.75] at (2.6,-0.15) {$,$};
        
            \node[fill=red] (r) at (6,-0.4) {};  
            \node (h2) at (6.259,0.566) {};
            \node (h3) at (6.966,-0.141) {};
        
            \draw[thick] (r)--(h2)--(h3)--(r);
        
            \node[fill=none, scale=1.75] at (5,0) {$H =$};
            \node[fill=none, scale=1.75] at (6,-0.75) {$r$};
        \end{tikzpicture}
    \end{center}
    Then, the rooted product yields the graph
    
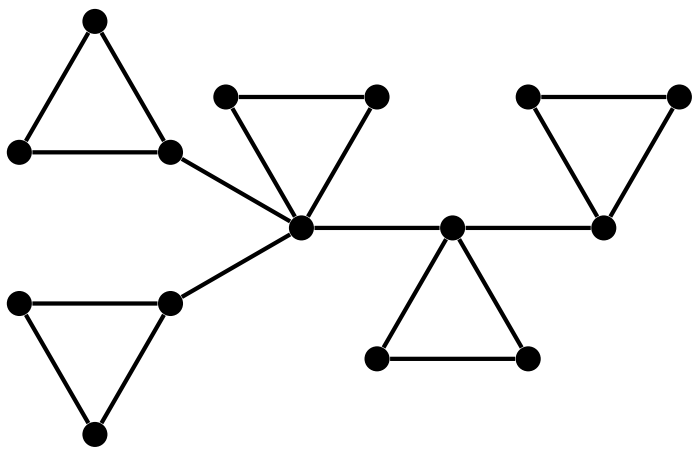
\begin{center}
        \begin{tikzpicture}[auto=center,every node/.style={circle, fill=black, scale=0.5}]
            \node (p1) at (0,0) {};          
            \node (p2) at (-0.866,0.5) {};   
            \node (p3) at (-0.866,-0.5) {};  
            \node (p4) at (1,0) {};          
            \node (p5) at (2,0) {};          
        
            \node (a1) at (-0.5,0.866) {};   
            \node (b1) at (0.5,0.866) {};
            \node (a2) at (-1.366,1.366) {}; 
            \node (b2) at (-1.866,0.5) {};
            \node (a3) at (-1.866,-0.5) {};  
            \node (b3) at (-1.366,-1.366) {};
            \node (a4) at (0.5,-0.866) {};   
            \node (b4) at (1.5,-0.866) {};
            \node (a5) at (1.5,0.866) {};    
            \node (b5) at (2.5,0.866) {};
        
            \draw[thick] (p1)--(p2);
            \draw[thick] (p1)--(p3);
            \draw[thick] (p1)--(p4)--(p5);
        
            \draw[thick] (p1)--(a1)--(b1)--(p1);
            \draw[thick] (p2)--(a2)--(b2)--(p2);
            \draw[thick] (p3)--(a3)--(b3)--(p3);
            \draw[thick] (p4)--(a4)--(b4)--(p4);
            \draw[thick] (p5)--(a5)--(b5)--(p5);
        
            \node[fill=none, scale=1.75] at (-3.3,0) {$G(H_r) =$};
        \end{tikzpicture}
    \end{center}
\end{example}

\begin{remark}
A couple of immediate remarks are in order:
    \begin{enumerate}
        \item $G(H)$ has $|V(G)| \cdot |V(H)|$ vertices.
        \item $G(H)$ has $|E(G)| + |V(G)| \cdot |E(H)|$ edges.
        \item When $G$ and $H$ are both trees, then $G(H)$ is a tree.
    \end{enumerate}
\end{remark}

We will now present three results about the rooted product of graphs. These results, though simple, are essential for our main result.  The first lemma is straightforward, so we omit the proof here. Recall that $G\setminus v$ is the graph obtained by deleting the vertex $v\in V(G)$ and all edges incident to it. 

\begin{lemma}\label{lemma:attachment vertex deletion}
    Let $G$ be a simple graph and $(H,r)$ a rooted connected simple graph. Let $g_0 \in V(G)$. Then, we have $G(H_r) \setminus (g_0,r) \cong (G\setminus g_0)(H_r) \; \sqcup (H \setminus r)$.\qed
\end{lemma}

\begin{example}
    We show an illustration of Lemma \ref{lemma:attachment vertex deletion}. Let $G$ and $(H,r)$ be the graphs from Example \ref{ex:graph attachment}, and let $g_0$ be the degree 3 vertex in $G$. The figure below shows $G(H_r) \setminus (g_0,r)$. The isomorphic copy of $H \setminus r$ is colored red, and $(G\setminus g_0) (H_r)$ is colored blue.
    \begin{center}
        \begin{tikzpicture}[auto=center,every node/.style={circle, fill=black, scale=0.5}]
            \node[blue] (p2) at (-0.866,0.5) {};   
            \node[blue] (p3) at (-0.866,-0.5) {};  
            \node[blue] (p4) at (1,0) {};          
            \node[blue] (p5) at (2,0) {};          
        
            \node[red] (a1) at (-0.5,0.866) {};   
            \node[red] (b1) at (0.5,0.866) {};
            \node[blue] (a2) at (-1.366,1.366) {}; 
            \node[blue] (b2) at (-1.866,0.5) {};
            \node[blue] (a3) at (-1.866,-0.5) {};  
            \node[blue] (b3) at (-1.366,-1.366) {};
            \node[blue] (a4) at (0.5,-0.866) {};   
            \node[blue] (b4) at (1.5,-0.866) {};
            \node[blue] (a5) at (1.5,0.866) {};    
            \node[blue] (b5) at (2.5,0.866) {};
        
            \draw[thick, blue] (p4)--(p5);
        
            \draw[thick,red] (a1)--(b1);
            \draw[thick,blue] (p2)--(a2)--(b2)--(p2);
            \draw[thick,blue] (p3)--(a3)--(b3)--(p3);
            \draw[thick,blue] (p4)--(a4)--(b4)--(p4);
            \draw[thick,blue] (p5)--(a5)--(b5)--(p5);
        
            \node[fill=none, scale=1.75] at (-3.6,0) {$G(H_r) \setminus (g_0,r) \cong $};
            
        \end{tikzpicture}
    \end{center}
\end{example}

In what follows, our goal is to prove that $G(H)$ is isomorphic to $G'(H)$ if and only if $G$ and $G'$ are isomorphic.  We begin with the following lemma. 

\begin{lemma}\label{lemma:attachment isomorphism structure}
    Let $G,G'$ be simple connected graphs with $n \geq 2$ vertices, and let $(H,r)$ be a rooted connected simple graph. Suppose there exists a graph isomorphism $\varphi:V(G(H_r)) \to V(G' (H_r))$. Then, for any $g \in V(G)$, we must have $\varphi(g,r) \in V(G') \times \{r\}$.
\end{lemma}
\begin{proof}
    Assume for the sake of contradiction that $\varphi(g,r) = (g',h)$ where $g' \in V(G')$ and $h\neq r$. By Lemma \ref{lemma:attachment vertex deletion}, we have
    \[
    G(H_r)\setminus (g,r) \cong (G \setminus g) (H_r) \; \sqcup (H\setminus r).
    \]
    Thus, the connected component of maximum order in $(G(H_r))\setminus (g,r)$ has at most $(n-1)\cdot |V(H)|$ vertices. On the other hand, since $\varphi(g,r) = (g',h)$ and $h \neq r$, it follows that $G'(H_r) \setminus \varphi(g,r)$ has a connected component with at least $(n-1)\cdot |V(H)| + 1$ vertices. Therefore, $G(H_r)\setminus (g,r) \ncong G' (H_r) \setminus \varphi(g,r)$, which is a contradiction since isomorphism is preserved under vertex deletion.
\end{proof}

The next result guarantees that isomorphism behaves well with the rooted product of graphs.

\begin{proposition}\label{prop:attachment isomorphism preserved}
    Let $G,G'$ be simple connected graphs with $n \geq 2$ vertices, and $(H,r)$ a rooted connected graph. Then, $G \cong G'$ if and only if $G(H_r) \cong G' (H_r)$.
\end{proposition}
\begin{proof}
    Suppose that $\varphi:V(G) \to V(G')$ is a graph isomorphism, and define $\tilde{\varphi}:V(G(H_r)) \to V(G'(H_r))$ by 
    \[
    \tilde{\varphi}(g,h) = (\varphi(g),h).
    \]
    Since $\varphi$ is a bijection, it is immediate from the definition of the map that $\tilde{\varphi}$ is a bijection. We now check edge-compatibility. Suppose that $(g_1,h_1)(g_2,h_2) \in E(G(H_r))$. There are two cases. If $g_1=g_2$, then $h_1h_2 \in E(H)$, and so
    \[
    \tilde{\varphi}(g_1,h_1)\tilde{\varphi}(g_2,h_2) = (\varphi(g_1),h_1)(\varphi(g_2),h_2) \in E(G'(H_r))
    \]
    since $\varphi(g_1)=\varphi(g_2)$. Otherwise, we must have $h_1=h_2=r$ and $g_1 \neq g_2$. Since $\varphi$ is an isomorphism between $G$ and $G'$, we have $\varphi(g_1)\varphi(g_2) \in E(G')$, and so $(\varphi(g_1),r)(\varphi(g_2),r) \in E(G' (H_r))$. The same argument with $\varphi^{-1}$ shows that the edges in $G' (H_r)$ are preserved under $\tilde{\varphi}^{-1}$.

    For the converse, assume that $\psi:V(G(H_r)) \to V(G' (H_r))$ is a graph isomorphism. Let $G_r$ denote the subgraph of $G(H_r)$ induced by the vertex set $\{(g,r) : g \in G\}$. Clearly, $G_r \cong G$ and $G'_r \cong G'$. By Lemma \ref{lemma:attachment isomorphism structure}, $\psi$ restricts to an isomorphism between $G_r$ and $G_r'$, and so $G \cong G'$.
\end{proof}

We are particularly interested in the rooted product when $G=T$ is a tree and $(H,r)$ 
is the $k$-star $St_k$ rooted at its central vertex. In the next lemma, we describe how the core interacts with the rooted product in this case.  For simplicity, when we write $T(St_k)$, we always mean the rooted product where $St_k$ is rooted at its central vertex.

\begin{lemma}\label{lemma:star attachment core}
    Let $T$ be a tree and $k \geq 2$. Then:
    \begin{enumerate}
        \item $K_{T (St_k)} \cong T$ as unweighted graphs,
        \item $w(v) = k$ for any $v \in V(K_{T(St_k)})$, and
        \item $\lc(T(St_k)) = (k^{|V(T)|})$ and so $T(St_k)$ is proper
    \end{enumerate}
\end{lemma}
\begin{proof}
    Let $V(T) = \{u_1,\ldots,u_n\}$ and $V(St_k)=\{r, \ell_1,\ldots,\ell_{k-1}\}$ where $r$ is the central vertex of $St_k$. For (1), note that a vertex in $T (St_k)$ is internal if and only if it is of the form $(u_j,r)$. Hence, $IV(T (St_k)) = \{(u_1,r),\ldots, (u_n,r)\}$. Then, $E(K_{T(St_k)}) = \{(u,r)(v,r) : uv \in E(T)\}$ and so it follows that $K_{T (St_k)} \cong T$. For (2), note that for any $j$, the leaves in the neighborhood of $(u_j,r)$ are precisely $\{(u_j,\ell_1),\ldots,(u_j,\ell_{k-1})\}$, so $w((u_j,r)) = k$. (3) follows from Remark \ref{rmk:lead partition and core weights} and the fact that $k \geq 2$.
\end{proof}

\begin{example}\label{ex:attachment to trees}
    Let $G=T$ be the tree from Example \ref{ex:graph attachment}, and $H=St_3$ rooted at its central vertex. Then, $T(St_3)$ is shown below, and it is readily seen that $\lc(T(St_3))=(3^5)$
    \begin{center}
        \begin{tikzpicture}[auto=center,every node/.style={circle, fill=black, scale=0.5}]
            \node (p1) at (0,0) {};          
            \node (p2) at (-0.866,0.5) {};   
            \node (p3) at (-0.866,-0.5) {};  
            \node (p4) at (1,0) {};          
            \node (p5) at (2,0) {};          
        
            \node (a1) at (-0.5,0.866) {};   
            \node (b1) at (0.5,0.866) {};
            \node (a2) at (-1.366,1.366) {}; 
            \node (b2) at (-1.866,0.5) {};
            \node (a3) at (-1.866,-0.5) {};  
            \node (b3) at (-1.366,-1.366) {};
            \node (a4) at (0.5,-0.866) {};   
            \node (b4) at (1.5,-0.866) {};
            \node (a5) at (1.5,0.866) {};    
            \node (b5) at (2.5,0.866) {};
        
            \draw[thick] (p1)--(p2);
            \draw[thick] (p1)--(p3);
            \draw[thick] (p1)--(p4)--(p5);
        
            \draw[thick] (p1)--(a1);
            \draw[thick] (b1)--(p1);
            \draw[thick] (p2)--(a2);
            \draw[thick] (b2)--(p2);
            \draw[thick] (p3)--(a3);
            \draw[thick] (b3)--(p3);
            \draw[thick] (p4)--(a4);
            \draw[thick] (b4)--(p4);
            \draw[thick] (p5)--(a5);
            \draw[thick] (b5)--(p5);
        
            \node[fill=none, scale=1.75] at (-3.3,0) {$T(St_3) =$};
        \end{tikzpicture}
    \end{center}
\end{example}

\subsection{Proper-Tree Reduction}

We now show our main result, that the Tree Isomorphism Conjecture for proper trees implies the full conjecture. With the same notation as above, for $(H,r)$ a rooted graph and each $i \in \mathbb{P}$, define
\[
R_i(\boldsymbol{x}) = \sum_{\substack{\kappa \in \pc(H) \\\kappa(r)=i}} \prod_{v \in V(H)} x_{\kappa(v)}.
\]
In other words, $R_i(\boldsymbol{x})$ is precisely the sum of terms of $\csf_H$ where the root $r$ has been assigned color $i$. Since the context is clear, we will simply denote $R_i(\boldsymbol{x})$ by $R_i$. In the following theorem we show that the chromatic symmetric function of $G(H_r)$ can be obtained from the chromatic symmetric function of $G$ via an evaluation of $x_i = R_i$. 

\begin{theorem}\label{thm:general attachment csf}
    For $G$ a simple graph, and $(H,r)$ a connected, rooted simple graph, we have
    \[
    \csf_{G(H_r)} = \csf_G(R_1,R_2,\ldots),
    \]
    where $\csf_G(R_1,R_2,\ldots)$ denotes the evaluation of $\csf_G$ at $x_i = R_i$.
\end{theorem}
\begin{proof}
   Let $V(G)=\{g_1,\ldots,g_n\}$ and $V(H) = \{h_1, \ldots, h_m\}$. Without loss of generality, let $h_1=r$ be the root of $H$. We let $H_{g_j}$ denote the copy of $H$ in $G(H_r)$ rooted at $g_j$ for $1\leq j \leq n$. More formally, we let $H_{g_j}$ be the subgraph of $G(H_r)$ induced by the vertices
   $
   \{(g_j,h_1),\ldots,(g_j,h_m)\}.
   $
   Given the following:
   \begin{enumerate}
       \item[(i)] A proper coloring $\alpha \in\pc(G)$, and
       \item[(ii)] $n$ proper colorings $\gamma_1,\ldots,\gamma_n \in\pc(H)$ such that $\gamma_i(r) = \alpha(g_i)$ for each $1 \leq i \leq n$,
   \end{enumerate}
    we can define a proper coloring $\kappa \in \pc(G(H_r))$ by letting $\kappa(g,r) = \alpha(g)$ and $\kappa(g_i,h) = \gamma_i(h)$ for $h\neq r$. This is of course well-defined since $\gamma_i(r)=\alpha(g_i)$. Similarly, given a proper coloring $\kappa\in\pc(G(H_r))$, it is clear that the restrictions $\alpha=\kappa|_{V(G) \times \{r\}}$ and $\gamma_i = \kappa|_{H_{g_i}}$ satisfy properties (i) and (ii) above.
    
    Thus, by this characterization of the proper colorings of $G(H_r)$, we obtain:
   \begin{align*}
       \csf_{G(H)} 
        &= \sum_{\kappa \in \pc(G(H_r))} \prod_{(g,h) \in V(G(H_r))} x_{\kappa(g,h)} \\
        &= \sum_{\kappa \in \pc(G(H_r))} (x_{\kappa(g_1,h_1)}\cdots x_{\kappa(g_1,h_m)}) \cdot \ldots \cdot (x_{\kappa(g_n,h_1)} \cdots x_{\kappa(g_n,h_m)}) \\
        &= \sum_{\kappa \in \pc(G(H_r))} \left(\prod_{i=1}^{m} x_{\kappa(g_1,h_i)}\right) \cdots  \left(\prod_{i=1}^m x_{\kappa(g_n,h_i)}\right) \\
        &= \sum_{\alpha\in\pc(G)}\; \sum_{\substack{\gamma_1\in\pc(H_{g_1}) \\ \gamma_1(r)=\alpha(g_1)}} \; \cdots \; \sum_{\substack{\gamma_n\in\pc(H_{g_n}) \\ \gamma_n(r)=\alpha(g_n)}}  \left(\prod_{i=1}^{m} x_{\gamma_1(g_1,h_i)}\right) \cdots  \left(\prod_{i=1}^m x_{\gamma_n(g_n,h_i)}\right) \\
        &= \sum_{\alpha\in\pc(G)} R_{\alpha(g_1)} \cdots R_{\alpha(g_n)} \\
        &= \csf_G(R_1,R_2,\ldots) \qedhere
   \end{align*} 
\end{proof}

Using Theorem \ref{thm:general attachment csf}, we are able to show that we can reduce the Tree Isomorphism Conjecture to proper trees. And that this result will imply it for all trees. Recall that if $T$ is a tree and $H_r$ is a rooted tree, then $T(H_r)$ is a tree. 

\begin{theorem}\label{thm:general tree reduction}
    Let $T,T'$ be trees, and $(H,r)$ a rooted tree. Then,
    \[
    \csf_T = \csf_{T'} \; \Longleftrightarrow  \; \csf_{T (H_r)} = \csf_{T' (H_r)}
    \]
\end{theorem}
\begin{proof}
    Let $V(T)=\{v_1,\ldots,v_n\}$, $V(T') = \{w_1,\ldots,w_n\}$, and $V(H) = \{h_1,\ldots,h_m\}$ where we can assume without loss of generality that $r=h_1$ is the root of $H$. The forward direction follows immediately from Theorem \ref{thm:general attachment csf}. We now prove the converse. Assume that $\csf_{T(H_r)} = \csf_{T' (H_r)}$, which by the same theorem implies
    \[
    \csf_T(R_1,R_2,\ldots) = \csf_{T'}(R_1,R_2,\ldots).
    \]
    For any partition $\lambda \vdash n$, let $
    c_\lambda(T) = \# \{ A \subseteq E(T) : \lambda(A) = \lambda \}$,
    and define $c_\lambda(T')$ analogously. Then, by Theorem \ref{thm:stanley-power-sum-expansion}, we have
    \begin{align*}
        \csf_T(R_1,R_2,\ldots) - \csf_{T'}(R_1,R_2,\ldots)
            = \sum_{\lambda \vdash n} (-1)^{n - \ell(\lambda)} (c_\lambda(T) - c_\lambda(T')) p_\lambda(R_1,R_2,\ldots) = 0.
    \end{align*}
    We show that $\{p_\lambda(R_1,R_2,\ldots) : \lambda \vdash n \} \subset \Lambda^{nm}$ is linearly independent, which will finish the proof. Suppose that there exist coefficients $a_\lambda \in \C$ not all zero such that
    \begin{equation}\label{eq:plambda independence}
        \sum_{\lambda\vdash n} a_\lambda p_\lambda(R_1,R_2,\ldots) = 0.
    \end{equation}
    Let $\lambda \vdash n$ be a longest partition such that $a_\lambda \neq 0$ and write $\lambda =(\lambda_1,\ldots,\lambda_\ell)$, so that
    \begin{equation}\label{eq:p evaluation}
    p_\lambda(R_1,R_2,\ldots) = \prod_{j=1}^\ell \left(\sum_{i=1}^\infty R_i^{\lambda_j}\right) 
    \end{equation}
    Note that expanding the right-hand side of \eqref{eq:p evaluation} into monomials in $\boldsymbol{x}$, each monomial corresponds to a sequence of independent choices: for each $1 \leq j \leq \ell$, choosing a color $i_j$ for the root $r$, and subsequently choosing $\lambda_j$ proper colorings of $H$ that assign the color $i_j$ to the root. Note that excluding the color of the root, each of the $\lambda_j$ proper colorings of $H$ uses at most $m-1$ distinct colors. Hence, any monomial in the right-hand side of \eqref{eq:p evaluation} has at most $\ell+ \sum_j \lambda_j (m-1) = \ell+n(m-1)$ distinct variables. Moreover, it is easy to see that choosing distinct colors $i_1,\ldots,i_\ell$ for the roots, and then choosing a new color for each of the other vertices yields a monomial with exactly $\ell+n(m-1)$ distinct variables, which is indexed by the partition $\tilde{\lambda}=(\lambda,1^{n(m-1)})$. Note that this argument also shows that the only monomials in $p_\lambda(R_1,R_2,\ldots)$ with exactly $\ell+n(m-1)$ distinct variables must be indexed by $\tilde{\lambda}$.

     Now let $f$ be any monomial whose exponents are indexed by $\tilde{\lambda}$ in the expansion of the LHS of \eqref{eq:plambda independence}, and let $c_f$ be its coefficient. Suppose $\mu \vdash n$ has $\ell(\mu) < \ell$ and $a_\mu \neq 0$. By the same argument as above, every monomial in $p_\mu(R_1,R_2,\ldots)$ has at most $\ell(\mu)+n(m-1)$ distinct variables. Since $\ell(\mu)<\ell$, then $p_\mu(R_1,R_2,\ldots)$ does not contribute to $c_f$ in the LHS of \eqref{eq:plambda independence}. Now suppose $\mu \vdash n$ has $a_\mu \neq 0$ and $\ell(\mu)= \ell$, and that $f$ appears in the monomial expansion of $p_\mu(R_1,R_2,\ldots)$. By the argument in the preceding paragraph, $f$ must be indexed by $\tilde{\mu}$, so we must have $\tilde{\mu}=\tilde{\lambda}$, and so $\mu = \lambda$.
    
    Therefore, it follows that $c_f$ is a nonzero multiple $a_\lambda$, which is a contradiction.
\end{proof}

Theorem \ref{thm:general tree reduction} yields a many reductions of Stanley's conjecture to several subclasses of trees. Our main result is that, in particular, the theorem above reduces Stanley's conjecture to proper trees. For any $k \geq 2$, let $\mathcal{P}_k$ denote the class of proper trees whose leaf components all have order $k$, which we refer to as \defn{$k$-weighted trees}.

\begin{corollary}[Proper-tree reduction]\label{cor:equal weight proper reduction}
    Fix any integer $k \geq 2$. Then, it suffices to check Stanley's conjecture for trees $T \in \mathcal{P}_k$. In particular, it suffices to prove Stanley's conjecture for proper trees.
\end{corollary}
\begin{proof}
    This follows from Lemma \ref{lemma:star attachment core}, Theorem \ref{thm:general attachment csf} and Theorem \ref{thm:general tree reduction} by letting $H=St_k$ and the root $r$ to be the central vertex of the star.
\end{proof}

\begin{figure}[h]
    \centering
    \begin{tikzpicture}[
    every node/.style={
        circle,
        draw,
        fill=white,
        minimum size=5mm,
        inner sep=1pt
    },
    thick
]
    \node (v1) at (0,0) {$k$};
    \node (v2) at (1.5,0) {$k$};
    \node (v3) at (3,0) {$k$};
    \node (v4) at (4.5,0) {$k$};
    \node (v5) at (6,0) {$k$};
    \node (v6) at (7.5,0) {$k$};

    \node (v7) at (1.5,1.2) {$k$};
    \node (v8) at (0.7,2.2) {$k$};
    \node (v9) at (2.3,2.2) {$k$};

    \node (v10) at (3,-1.2) {$k$};
    \node (v11) at (2.2,-2.2) {$k$};
    \node (v12) at (3.8,-2.2) {$k$};

    \node (v13) at (4.5,1.2) {$k$};
    \node (v14) at (5.5,2.2) {$k$};

    \node (v15) at (6,-1.2) {$k$};

    \draw (v1)--(v2)--(v3)--(v4)--(v5)--(v6);
    \draw (v2)--(v7);
    \draw (v8)--(v7)--(v9);
    \draw (v3)--(v10);
    \draw (v11)--(v10)--(v12);
    \draw (v4)--(v13)--(v14);
    \draw (v5)--(v15);
\end{tikzpicture}
    \caption{A general tree in $\mathcal{P}_k$, represented by its weighted core.}
    \label{fig:pk tree}
\end{figure}
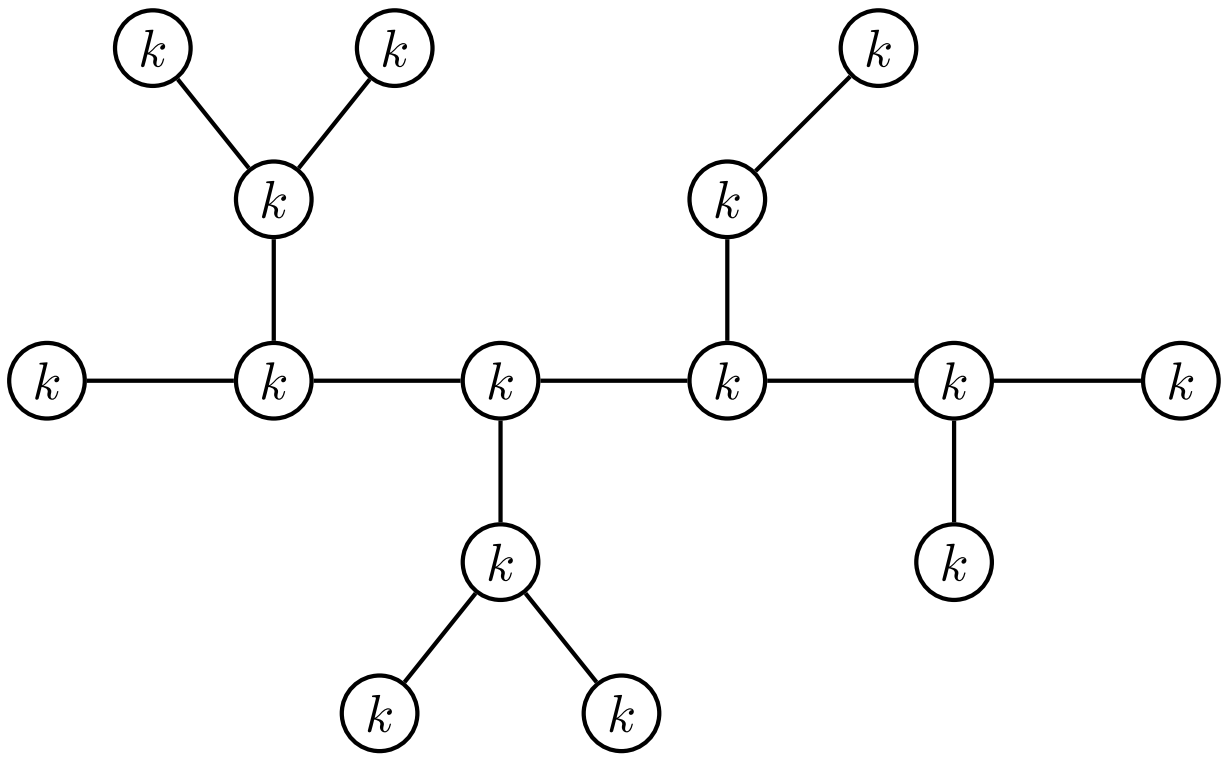

In Section \ref{sec:k weighted trees}, we will further explore the structure of certain coefficients in the chromatic symmetric function of trees that belong to $\mathcal{P}_k$, making the relevance of Corollary \ref{cor:equal weight proper reduction} even more clear. For a general example of a tree in $\mathcal{P}_k$, see Figure \ref{fig:pk tree}.

\section{$k$-weighted trees}\label{sec:k weighted trees}

In this section, we study the chromatic symmetric function of trees in $\mathcal{P}_k$, where $k \geq 2$. Recall that, by Corollary \ref{cor:equal weight proper reduction}, the Tree Isomorphism Conjecture is equivalent to verifying it in the class $\mathcal{P}_k$. We first prove that in $\mathcal{P}_k$ the information about the tree carried by some coefficients can be obtained from the subset of terms that do not contain 1s. That is, the coefficients of the chromatic symmetric function in $\mathcal{P}_k$ carry redundant information. Then, we show a correspondence between the 1-free layer of coefficients in the $\stfrak_\lambda$-expansion of $\csft$ and the signed $p_\lambda$-expansion of their unweighted core.

\subsection{1-free coefficients} Here, we show that if $T \in \mathcal{P}_k$, then the coefficients $c_\lambda$ in the star-basis expansion of $\csft$ with $1 \notin \lambda$ in fact determine $\csft$ completely. Loosely stated, this guarantees that in $\mathcal{P}_k$, deletions and leaf-contractions carry all the information about the tree $T$ in $\csft$.

Before stating and proving the theorem, we make a brief remark. For $T \in \mathcal{P}_k$, every leaf component in $T$ has order $k \geq 2$. Hence, any sequence of operations in a DNC tree for $T$ obtained by applying only deletions and leaf-contractions will result in a star forest $F$ where $\lambda(F)$ contains only parts that are multiples of $k$. In particular, $1 \notin \lambda(F)$.

\begin{theorem}\label{thm:1-free coefficients}
    Let $T \in \mathcal{P}_k$ with $k \geq 2$ and let $\csft = \sum_\mu c_\mu \stfrak_\mu$ be the star-basis expansion of its chromatic symmetric function. Then, for any $\mu$ with $c_\mu \neq 0$, we must have $\mu=(\mu_1, \ldots, \mu_\ell, k^{m_k}, 1^{m_1} )$ where $\mu_1 \geq \cdots \geq \mu_\ell > k$. Moreover, if $m_1 \geq 1$, then
    \begin{equation}\label{eq:coefficients with 1}
        c_\mu = (-1)^{m_1}\sum_{\mathbf{i}=(i_1,\ldots,i_\ell)} c_{\mu^\mathbf{i}} \cdot \prod_{\substack{p>k \\ m_p(\mu^\mathbf{i})>0}} \binom{m_p(\mu^{\mathbf{i}})}{\eta_{p,0},\eta_{p,1},\ldots, \eta_{p, m_1}} \prod_{d:\eta_{p,d}\neq 0} \binom{p/k-1}{d}^{\eta_{p,d}}
    \end{equation}
    where $\mu^\mathbf{i} = \sort(\mu_1+i_1, \ldots, \mu_\ell+i_\ell, k^{m_k})$, $\eta_{p,d}=\#\{j \in [\ell] : \mu_j+i_j = p, i_j = d\}$
    and the sum runs over $(i_1,\ldots,i_\ell) \in \Z_{\geq 0}^\ell$ such that
    \begin{itemize}
        \item[(i)] $\sum_j i_j = m_1$, 
        \item[(ii)] Each $\mu_j + i_j$ is a multiple of $k$, 
        \item[(iii)] and $i_j \geq i_{j+1}$ whenever $\mu_j=\mu_{j+1}$.
    \end{itemize}     
\end{theorem}
\begin{proof}
    Fix a DNC tree $\mathcal{T}$ for $T$, and let $S_\mu$ be the set of DNC sequences in $\mathcal{T}$ that end in a star forest $F$ indexed by $\mu$. Then, by Lemma \ref{lem:coeff-paths-seq}, $c_\mu = (-1)^{m_1}|S_\mu|$. We now show that $|S_\mu|$ is  the sign-free quantity in Equation \eqref{eq:coefficients with 1}.

    Since $T$ is proper, any sequence in $S_\mu$ has length $|I(T)|$. Moreover, since each leaf component in $T$ has order $k$, any such sequence satisfies the following:
    \begin{itemize}
        \item The number of $L$s is $\ell+m_k-1$,
        \item The number of $M$s is $m_1$, and
        \item The number of $R$s is $|I(T)|-m_1-(\ell+m_k-1)$.
    \end{itemize}
    Now fix a sequence $s \in S_\mu$, and let $e_1,\ldots,e_{\ell+m_k-1} \in I(T)$ be the internal edges in $T$ that correspond to an $L$ in $s$. Then, $T_s = T \setminus \{e_1,\ldots,e_{\ell+m_k-1}\}$ has $\ell+m_k$ connected components. Any connected component having $t \geq 1$ internal edges induces a star of order at least $(t+1)k-t > k$ in $F$. Therefore, since $\lambda(F)$ has $m_k$ parts of size $k$, exactly $m_k$ connected components in $T_s$ are stars of order $k$. Hence, the $m_1$ dot-contractions must be applied to internal edges that are contained in the remaining $\ell$ connected components of $T_s$. Order these $\ell$ connected components by decreasing order of the stars that they induce on $F$, and by decreasing number of dot-contractions whenever two components induce a star of the same order in $F$. Denote these components by $C_1,\ldots,C_\ell$, and let $i_j$ be the number of dot-contractions that are applied in $C_j$. Clearly we have $\sum_j i_j = m_1$, and $|C_j|-i_j = \mu_j$, so each $\mu_j+i_j$ is a multiple of $k$. Moreover, by the ordering we imposed, $i_j \geq i_{j+1}$ whenever $\mu_j = \mu_{j+1}$.

    Thus, for any sequence $s \in S_\mu$, we have a list $\mathbf{i}_s=(i_1,\ldots,i_\ell)$ satisfying the three conditions in the statement of the theorem. Then, the sets $S_\mu(\mathbf{i}) = \{s \in S_\mu : \mathbf{i}_s = \mathbf{i}\}$ partition $S_\mu$, so
    \[
    |S_\mu|= \sum_{\mathbf{i}} |S_{\mu}(\mathbf{i})|
    \]
    where the sum runs over the lists satisfying conditions (i), (ii) and (iii). 
    
    Now fix a list $\mathbf{i}$ satisfying the three conditions above. We count $|S_\mu(\mathbf{i})|$ as follows: (1) choose the set of edges where an $L$ is applied, (2) assign a number of dot-contractions to each resulting connected component, and (3) choose the edges on which these dot-contractions occur within each connected component.
    \begin{enumerate}
        \item \emph{Choose the $L$-edges.}
        The set $E \subseteq I(T)$ of deleted edges must satisfy
        $\lambda(T\setminus E)=\mu^{\mathbf{i}}$. Thus, the choice of $L$-edges ranges over the set
        \[
        \mathcal{D}_{\mu^{\mathbf{i}}} = \{E\subseteq I(T):\lambda(T\setminus E)=\mu^{\mathbf{i}}\},
        \]
        and it is straightforward to see that $|\mathcal{D}_{\mu^{\mathbf{i}}}| = c_{\mu^{\mathbf{i}}}$.
        
        \item \emph{Assign number of dot-contractions to connected components.}
        Fix $E\in\mathcal D_{\mu^{\mathbf{i}}}$ from step (1). For each $p>k$, the forest
        $T\setminus E$ has $m_p(\mu^{\mathbf i})$ components of order $p$.
        To obtain the list of dot-contraction counts $\mathbf{i}$, we have to apply $d$ dot-contractions in exactly $\eta_{p,d}$ of
        these components, for each $0\leq d\leq m_1$. Since $\sum_{d=0}^{m_1}\eta_{p,d}=m_p(\mu^{\mathbf i})$,
        the number of ways to assign the number of dot-contractions given to the components of order $p$ is
        \[
            \binom{m_p(\mu^{\mathbf i})}
            {\eta_{p,0},\eta_{p,1},\ldots,\eta_{p,m_1}}.
        \]
        Since we are assigning these dot-contractions independently for each component order, the total number of ways in which we can assign dot-contractions to the components is
        \[
            \prod_{\substack{p>k\\m_p(\mu^{\mathbf i})>0}}
            \binom{m_p(\mu^{\mathbf i})}
            {\eta_{p,0},\eta_{p,1},\ldots,\eta_{p,m_1}}.
        \] 
        \item \emph{Choose the dot-contracted edges.}
        Fix an assignment from step (2). A component $C$ of order $p$ in $T \setminus E$ consists of $p/k$ leaf components of order $k$, joined by $p/k-1$ internal edges. If $C$ has been assigned $d$ dot-contractions, we choose the dot-contracted edges in $\binom{p/k-1}{d}$ ways, and apply leaf-contractions (i.e. $R$s in the sequence) to the remaining internal edges in $C$.
    
        There are $\eta_{p,d}$ components of order $p$ to which we have to apply $d$ dot-contractions. Therefore, the number of choices within all components is\footnote{Note that in this product formula we are using the convention that $\binom{a}{d}=0$ when $d>a$}
        \[
            \prod_{\substack{p>k\\m_p(\mu^{\mathbf i})>0}}
            \prod_{\substack{0\leq d\leq m_1\\\eta_{p,d}>0}}
            \binom{p/k-1}{d}^{\eta_{p,d}}.
        \]
        
    \end{enumerate}

    Now note that these three steps determine the operation applied to every
    internal edge in $T$. More specifically, the edges in $E$ receive $L$, the chosen
    dot-contraction edges receive $M$, and all remaining internal
    edges receive $R$. Moreover, following these choices $\mathcal{T}$
    determines a \emph{unique} sequence in $S_\mu(\mathbf{i})$. Conversely, every sequence in $S_\mu(\mathbf{i})$ uniquely determines its deletion set, the contraction count assigned to
    each component, and the dot-contraction edges within each
    component. Thus, our succession of choices counts every sequence exactly
    once.
    
    Lastly, observe that the number of choices in steps (2) and (3) depends only
    on $\mathbf{i}$ and $\mu^{\mathbf{i}}$, and is thus the same for every
    $E\in\mathcal D_{\mu^{\mathbf{i}}}$. Hence, we obtain
    \[
        |S_\mu(\mathbf{i})|
        =
        c_{\mu^{\mathbf{i}}}
        \prod_{\substack{p>k\\m_p(\mu^{\mathbf{i}})>0}}
            \binom{m_p(\mu^{\mathbf{i}})}
            {\eta_{p,0},\eta_{p,1},\ldots,\eta_{p,m_1}}
            \prod_{\substack{0\leq d\leq m_1\\\eta_{p,d}>0}}
            \binom{p/k-1}{d}^{\eta_{p,d}}.
    \]
\end{proof}

We now give an example to illustrate the computation using Equation \eqref{eq:coefficients with 1}.

\begin{example}
    We give an example of Theorem \ref{thm:1-free coefficients} when $k=2$. We will first show a direct computation using the expression in \eqref{eq:coefficients with 1}, and then comment on the intuition behind it.
    
    First, consider the following chromatic symmetric function 
        \[
        \begin{aligned}
        X_T ={}&
        \stfrak_{(12)} - 5\stfrak_{(11,1)} + 3\stfrak_{(10,2)}
        + 10\stfrak_{(10,1^2)} - 12\stfrak_{(9,2,1)}
        - 10\stfrak_{(9,1^3)} \\
        &+ \stfrak_{(8,4)} - \stfrak_{(8,3,1)} + 4\stfrak_{(8,2^2)}
        + 18\stfrak_{(8,2,1^2)} + 5\stfrak_{(8,1^4)}
        - 3\stfrak_{(7,4,1)} \\
        &+ 3\stfrak_{(7,3,1^2)} - 12\stfrak_{(7,2^2,1)}
        - 12\stfrak_{(7,2,1^3)} - \stfrak_{(7,1^5)}
        + \stfrak_{(6^2)} - 4\stfrak_{(6,5,1)} \\
        &+ 6\stfrak_{(6,4,2)} + 5\stfrak_{(6,4,1^2)}
        - 6\stfrak_{(6,3,2,1)} - 3\stfrak_{(6,3,1^3)}
        + 5\stfrak_{(6,2^3)} + 12\stfrak_{(6,2^2,1^2)} \\
        &+ 3\stfrak_{(6,2,1^4)} + 4\stfrak_{(5^2,1^2)}
        - 12\stfrak_{(5,4,2,1)} \textcolor{darkred}{- 5\stfrak_{(5,4,1^3)}}
        + 12\stfrak_{(5,3,2,1^2)} + \stfrak_{(5,3,1^4)} \\
        &- 10\stfrak_{(5,2^3,1)} - 4\stfrak_{(5,2^2,1^3)}
        + 5\stfrak_{(4^2,2^2)} + 6\stfrak_{(4^2,2,1^2)}
        + \stfrak_{(4^2,1^4)} - 10\stfrak_{(4,3,2^2,1)} \\
        &- 6\stfrak_{(4,3,2,1^3)} + 5\stfrak_{(4,2^4)}
        + 5\stfrak_{(4,2^3,1^2)} + 5\stfrak_{(3^2,2^2,1^2)}
        - 5\stfrak_{(3,2^4,1)} + \stfrak_{(2^6)}.
        \end{aligned}
        \]
        Since $\lead(\csft) = (2^6)$, then $T \in \mathcal{P}_2$. We verify Equation \eqref{eq:coefficients with 1} for $\mu=(5,4,1^3)$, which has been colored above. The only vectors satisfying the conditions from the theorem are $\mathbf{i}_1=(1,2)$ and $\mathbf{i}_2=(3,0)$.
        For the former, we have $\mu^{\mathbf{i}_1}=(6,6)$ and $c_{\mu^{\mathbf{i}_1}}=1$. The only $p>2$ with positive multiplicity in $\mu^{\mathbf{i}_1}$ is $p=6$, and in this case we have $\eta_{6,1}=1=\eta_{6,2}$ and $\eta_{6,0}=\eta_{6,3}=0$. To obtain a component of order 6, each connected component must have three leaf components of order 2, so it contains $2=6/2-1$ internal edges. Then, for $\mathbf{i}_1=(1,2)$, the summation term is
        \[
        c_{(6,6)} \cdot \binom{2}{1,1} \cdot \binom{2}{1} \cdot \binom{2}{2} = 4.
        \]
        
        For the latter choice of vector, we have $\mu^{\mathbf{i}_2}=(8,4)$ with $c_{(8,4)}=1$. In this case, there are two choices for $p$, namely $p=8$ and $p=4$. For $p=8$, we have $m_8(\mu^{\mathbf{i}_2})=1$ and $\eta_{8,3}=1$, while $\eta_{8,d}=0$ for any $d \neq 3$. To obtain a component of order 8, we need $3=8/2-1$ internal edges connecting leaf components of order 2. For $p=4$, we have $m_4(\mu^{\mathbf{i}_2}) = 1$ and $\eta_{4,0}=1$. Similarly as above, only one internal edge is required to obtain a component of order 4 via leaf components of order 2. Then, for $\mathbf{i}_2=(3,0)$, the summation term is
        \[
        c_{(8,4)} \cdot \left[\binom{1}{0,0,0,1} \cdot \binom{3}{3} \right] \cdot \left[\binom{1}{1} \cdot \binom{1}{0}\right] = 1.
        \]
        Then, we have obtained $c_{(5,4,1^3)}=(-1)^3 \cdot (4+1)=-5$, as expected.

        We now comment on the intuition behind this computation. It can be readily checked that $\csft$ above is the chromatic symmetric function of the tree below
        
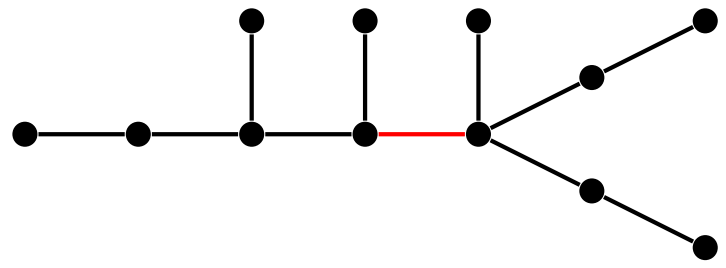
\begin{center}
            \begin{tikzpicture}[auto=center, every node/.style={circle, fill=black, scale=0.5}, thick, scale=0.75]
            \node (A11) at (-3,0) {};
            \node (A5)  at (-2,0) {};
            \node (A4)  at (-1,0) {};
            \node (A0)  at (0,0) {};
            \node (A1)  at (1,0) {};
        
            \node (A10) at (-1,1) {};
            \node (A6)  at (0,1) {};
            \node (A7)  at (1,1) {};
        
            \node (A2) at (2,0.5) {};
            \node (A8) at (3,1) {};
            \node (A3) at (2,-0.5) {};
            \node (A9) at (3,-1) {};
        
            \draw (A11)--(A5)--(A4);
            \draw (A4)--(A0);
            \draw[red] (A0)--(A1);
            \draw (A4)--(A10);
            \draw (A0)--(A6);
            \draw (A1)--(A7);
            \draw (A1)--(A2)--(A8);
            \draw (A1)--(A3)--(A9);
        \end{tikzpicture}
        \end{center}

        In the case $\mathbf{i}_1=(1,2)$, we obtained $c_{\mu^{\mathbf{i}_1}}=1$, which indicates that there is only one way of choosing a set of $L$-edges so that the resulting forest is indexed by $(6,6)$. This is clear since the only such edge is the edge $e$, which is colored red above. $T\setminus e$ thus has two connected components of order $6$. One has to apply 1 dot-contraction in one of said components, and 2 dot-contractions in the other to obtain the partition $(5,4,1^3)$. This illustrates that the meaning of the vectors $\mathbf{i}$ is just an assignment of dot-contraction counts to the connected components in $T \setminus E$, as detailed in Step (2) of the proof of Theorem \ref{thm:1-free coefficients}. Once the $L$-edges and the dot-contraction counts have been chosen, it only remains to choose the specific edges that are dot-contracted, and this can be done in $\binom{p/k-1}{d}$ ways as $d$ ranges over the dot-contraction counts in $\mathbf{i}$.
\end{example}

Lastly, we show an example that illustrates that a statement such as Theorem \ref{thm:1-free coefficients} does not hold outside of the class $\mathcal{P}_k$. More specifically, outside of $\mathcal{P}_k$, the coefficients of terms containing 1s contain strictly more information than those without 1s, even if we restrict ourselves to studying \emph{proper} trees outside of $\mathcal{P}_k$. The example suggests that approaching the isomorphism conjecture by proving it only on $\mathcal{P}_k$ is a worthwhile endeavor.

\begin{example}
    Consider the following trees.
    \vskip .1in
    \begin{center}
        \begin{tikzpicture}[
            auto=center,
            every node/.style={circle, fill=black, scale=0.5},
            thick
        ]
            \begin{scope}
                \node (g1) at (0,0) {};
                \node (g2) at (1.5,0) {};
                \node (g3) at (3,0) {};
                \node (g4) at (4.5,0) {};
                \node (g5) at (6,0) {};
        
                \node (g6) at (0,0.6) {};
        
                \node (g7) at (1.2,0.5) {};
                \node (g8) at (1.5,0.6) {};
                \node (g9) at (1.8,0.5) {};
        
                \node (g10) at (3,0.6) {};
        
                \node (g11) at (4,0.3) {};
                \node (g12) at (4.2,0.5) {};
                \node (g13) at (4.5,0.6) {};
                \node (g14) at (4.8,0.5) {};
                \node (g15) at (5,0.3) {};
        
                \node (g16) at (5.7,0.5) {};
                \node (g17) at (6,0.6) {};
                \node (g18) at (6.3,0.5) {};
        
                \draw (g1)--(g2)--(g3)--(g4)--(g5);
                \draw (g1)--(g6);
        
                \foreach \i in {7,8,9}
                    \draw (g2)--(g\i);
        
                \draw (g3)--(g10);
        
                \foreach \i in {11,12,13,14,15}
                    \draw (g4)--(g\i);
        
                \foreach \i in {16,17,18}
                    \draw (g5)--(g\i);
        
                \node[fill=none, scale=2] at (3.25,-0.6) {$T_1$};
            \end{scope}
        
            \begin{scope}[xshift=8cm]
                \node (h1) at (0,0) {};
                \node (h2) at (1.5,0) {};
                \node (h3) at (3,0) {};
                \node (h4) at (4.5,0) {};
                \node (h5) at (6,0) {};
        
                \node (h6) at (0,0.6) {};
        
                \node (h7) at (1,0.3) {};
                \node (h8) at (1.2,0.5) {};
                \node (h9) at (1.5,0.6) {};
                \node (h10) at (1.8,0.5) {};
                \node (h11) at (2,0.3) {};
        
                \node (h12) at (2.7,0.5) {};
                \node (h13) at (3,0.6) {};
                \node (h14) at (3.3,0.5) {};
        
                \node (h15) at (4.5,0.6) {};
        
                \node (h16) at (5.7,0.5) {};
                \node (h17) at (6,0.6) {};
                \node (h18) at (6.3,0.5) {};
        
                \draw (h1)--(h2)--(h3)--(h4)--(h5);
                \draw (h1)--(h6);
        
                \foreach \i in {7,8,9,10,11}
                    \draw (h2)--(h\i);
        
                \foreach \i in {12,13,14}
                    \draw (h3)--(h\i);
        
                \draw (h4)--(h15);
        
                \foreach \i in {16,17,18}
                    \draw (h5)--(h\i);
        
                \node[fill=none, scale=2] at (3.25,-0.6) {$T_2$};
            \end{scope}
        \end{tikzpicture}
    \end{center}

    Observe that both $T_1$ and $T_2$ are proper, have the same leaf component partition, and both are in fact caterpillars (i.e. the subgraph induced by the internal vertices is a path). However, one can compute directly that $\csf_{T_1}$ and $\csf_{T_2}$ agree exactly on all the terms whose indexing partitions do not contain 1s. Explicitly, these terms are
    \[
        \begin{aligned}
            &\stfrak_{(18)}
            +\stfrak_{(16,2)}
            +\stfrak_{(14,4)}
            +\stfrak_{(12,6)}
            +2\stfrak_{(12,4,2)}
            +\stfrak_{(10,8)} \\
            &\quad
            +2\stfrak_{(10,6,2)}
            +\stfrak_{(10,4,2,2)}
            +2\stfrak_{(8,6,4)}
            +\stfrak_{(8,4,4,2)}
            +2\stfrak_{(6,6,4,2)}
            +\stfrak_{(6,4,4,2,2)}.
        \end{aligned}
    \]
    Thus, $\csf_{T_1}$ and $\csf_{T_2}$ differ only on terms containing 1s. For instance, $[\stfrak_{(10,7,1)}]\csf_{T_1}=-2$, but $[\stfrak_{(10,7,1)}]\csf_{T_2}=-1$. Thus, outside of $\mathcal{P}_k$, the terms whose indexing partition contains 1 carry strictly more information than the 1-free terms.
\end{example}

\subsection{A $p_\lambda \leftrightarrow \stfrak_\lambda$ correspondence} Here, we present an alternative method to prove that, for $k$-weighted trees, the 1-free coefficients in the $\stfrak$-expansion encode all the information about $\csft$. Namely, we show a direct correspondence between the 1-free $\stfrak$-terms of $\csft$ and the power-sum expansion of $\csf_K$, where $K$ is the unweighted core of $T$.

Observe that any tree $T \in \mathcal{P}_k$ is of the form $T=K(St_k)$, where $K = T|_{IV(T)}$ is the unweighted core of $T$. Letting $V(St_k) = \{r,\ell_1,\ldots,\ell_{k-1}\}$ and   
\[
R_i = \sum_{\substack{\kappa \in \pc(St_k) \\ \kappa(r)=i}} x_{\kappa(r)} x_{\kappa(\ell_1)}\cdots x_{\kappa(\ell_{k-1})}~,
\]
Theorem \ref{thm:general attachment csf} then shows that $\csft=\csf_K(R_1,R_2,\ldots)$, so $\csft$ is fully determined by $\csf_K$ whenever $T$ is $k$-weighted.

Now since $T$ is proper, Proposition \ref{prop:csf in terms of core} implies that the $1$-free layer of the star-basis expansion of $\csft$, which we denote by $[\csft]_{\stfrak_1=0}$, satisfies
\[
[\csft]_{\stfrak_1=0} = \sum_{A \subseteq E(K)} \prod_{B \in c(A)} \stfrak_{w(B)} = \sum_{A\subseteq E(K)} \prod_{B \in c(A)} \stfrak_{k|B|} = \sum_{A \subseteq E(K)} \stfrak_{k \cdot \lambda(A)},
\]
where $k \cdot \lambda$ denotes the partition obtained by multiplying each part of $\lambda$ by $k$. Moreover, Theorem \ref{thm:stanley-power-sum-expansion} implies
\[
\csf_K = \sum_{A \subseteq E(K)} (-1)^{|V(K)|-\ell(\lambda(A))}p_{\lambda(A)}.
\]
Then, the 1-free layer of the $\stfrak$-expansion of $\csft$ completely determines $\csf_K$ which in turn completely determines $\csft$. Therefore, the terms $c_\lambda\stfrak_\lambda$ with $1 \in \lambda$ carry redundant information.

\begin{example}
    Let $T$ be the five-vertex tree from Example \ref{ex:graph attachment}. In the power-sum basis, we have
    \[
    \csft = p_{(1^5)} -4p_{(2,1^3)}+2p_{(2^2,1)}+4p_{(3,1^2)}- p_{(3,2)} - 3p_{(4,1)}+p_{(5)}.
    \]
    Now consider $T(St_3)$, which from the computation in Example \ref{ex:csf in terms of core} has
    \[
    [\csf_{T(St_3)}]_{\stfrak_1=0} = \stfrak_{(3^5)}
            +4\stfrak_{(6,3^3)}
            +2\stfrak_{(6^2,3)}
            +4\stfrak_{(9,3^2)}
            +\stfrak_{(9,6)}
            +3\stfrak_{(12,3)}
            +\stfrak_{(15)}.
    \]
    Then, it is clear that $[\csf_{T(St_3)}]_{\stfrak_1=0}$ is obtained by identifying $p_\lambda \leftrightarrow \stfrak_{3\cdot\lambda}$ and forgetting the signs in the power-sum expansion.
\end{example}

\printbibliography

\end{document}